\documentclass[final]{siamltex}

\usepackage{graphicx}
\usepackage{amsmath}
\usepackage{amssymb}
\usepackage{amsfonts}
\usepackage{amsxtra}
\usepackage{epsfig}
\usepackage{verbatim}
\usepackage{alltt}
\usepackage{bm}
\usepackage{caption}
\usepackage{subcaption}
\usepackage{algorithmicx}
\usepackage[boxed]{algorithm}
\usepackage{mathtools}
\usepackage{algpseudocode}
\numberwithin{equation}{section}
\usepackage{todonotes}
\usepackage{bm}
\usepackage{tikz}
\usepackage{pgfplots}
\usetikzlibrary{plotmarks}
\usepackage{hyperref}
\usepackage{framed}
\usepackage{minibox}

\floatname{algorithm}{Alg.}

\algtext*{EndProcedure}

\newcommand{\vecname}[1]{\ensuremath{\mathbf{#1}}}

\newcommand{\vx}{\vecname{x}} 
\newcommand{\vb}{\vecname{b}} 
\newcommand{\vy}{\vecname{y}} 
\newcommand{\vs}{\vecname{s}} 
\newcommand{\vd}{\vecname{d}} 
\newcommand{\vf}{f} 
\newcommand{\vr}{\vecname{r}} 
\newcommand{\vc}{\vecname{c}} 
\newcommand{\vF}{\vecname{F}} 
\newcommand{\vJ}{\vecname{J}} 
\newcommand{\vK}{\vecname{K}} 
\newcommand{\vM}{\vecname{M}} 
\newcommand{\vN}{\vecname{N}} 
\newcommand{\vA}{\vecname{A}} 
\newcommand{\vP}{\vecname{P}} 
\newcommand{\vQ}{\vecname{Q}} 
\newcommand{\vI}{\vecname{I}} 
\newcommand{\inj}{\vecname{\hat{R}}}
\newcommand{\rsrct}{\vecname{R}}
\newcommand{\intrp}{\vecname{P}}

\newcommand{\vu}{\vecname{u}} 

\newcommand {\pc}{\ensuremath{-}}
\newcommand {\lp}{\ensuremath{\pc_{L}}}
\newcommand {\rp}{\ensuremath{\pc_{R}}}
\newcommand {\lin}{\ensuremath{\backslash}}

\newcommand{\solvername}[1]{\ensuremath{\begingroup\text{#1}\endgroup}}

\newcommand {\NRICH}{\solvername{NRICH}}
\newcommand {\NGMRES}{\solvername{NGMRES}}
\newcommand {\ANDERSON}{\solvername{ANDERSON}}
\newcommand {\NEWT}{\solvername{NEWT}}
\newcommand {\NCG}{\solvername{NCG}}
\newcommand {\FAS}{\solvername{FAS}}
\newcommand {\QN}{\solvername{QN}}
\newcommand {\NASM}{\solvername{NASM}}
\newcommand {\RAS}{\solvername{RAS}}
\newcommand {\ASM}{\solvername{ASM}}
\newcommand {\MG}{\solvername{MG}}
\newcommand {\GSN}{\solvername{GSN}}
\newcommand {\GS}{\solvername{GS}}

\newcommand {\SOR}{\solvername{SOR}}
\newcommand {\ASPIN}{\solvername{ASPIN}}

\newcommand {\CG}{\solvername{CG}}
\newcommand {\GMRES}{\solvername{GMRES}}
\newcommand {\krylov}{\solvername{K}}

\newcommand {\LU}{\solvername{LU}}

\newcommand {\NK}{\NEWT\lin\krylov}

\newcommand {\LT}{\ensuremath{\text{L2}}}
\newcommand {\CP}{\ensuremath{\text{CP}}}
\newcommand {\BT}{\ensuremath{\text{BT}}}

\newcommand{\eqnref}[1]{(\ref{#1})}  
\newcommand{\figref}[1]{Fig.~\ref{#1}}  
\newcommand{\agmref}[1]{Alg.~\ref{#1}} 
\newcommand{\tabref}[1]{Table~\ref{#1}} 
\newcommand{\secref}[1]{Section~\ref{#1}} 

\newcommand {\composition}{composition }
\newcommand {\compositions}{compositions }

\begin{document}

%
%

\title{\bf Composing Scalable Nonlinear Algebraic Solvers}

\renewcommand{\thefootnote}{\arabic{footnote}}

\author{Peter R. Brune\thanks{prbrune@mcs.anl.gov and bsmith@mcs.anl.gov, Mathematics and Computer Science Division, Argonne National Laboratory, 9700 S. Cass Ave., Argonne, IL 60439} \and
Matthew G. Knepley\thanks{knepley@rice.edu, Department of Computational and Applied Mathematics, Rice University, Duncan Hall, 6100 S. Main St., Houston, TX 77005} \and
Barry F. Smith\footnotemark[1] \and
Xuemin Tu\thanks{xtu@math.ku.edu, Department of Mathematics, University of Kansas, 1460 Jayhawk Blvd., Lawrence, KS 66045}}

\maketitle
\begin{abstract}
  Most efficient linear solvers use composable algorithmic components, with the most common model being the combination
  of a Krylov accelerator and one or more preconditioners. A similar set of concepts may be used for nonlinear algebraic
  systems, where nonlinear composition of different nonlinear solvers may significantly improve the time to solution.
  We describe the basic concepts of nonlinear \composition and preconditioning and present a number of solvers
  applicable to nonlinear partial differential equations. We have developed a software framework in order to easily
  explore the possible combinations of solvers.  We show that the performance gains from using composed solvers can be
  substantial compared with gains from standard Newton-Krylov methods.
\end{abstract}
{\center
}
\section{Introduction}
\label{sec:intro}
Large-scale algebraic solvers for nonlinear partial differential equations (PDEs) are an essential component of modern
simulations.  Newton-Krylov methods \cite{DemboIN1982} have well-earned dominance. They are generally robust and may be
built from preexisting linear solvers and preconditioners, including fast multilevel preconditioners such as multigrid
\cite{bramble1,brandt1984,BriggsHensonMcCormick00,trottenberg2001multigrid,wessling1992} and domain decomposition
methods \cite{qv99,1sbg,st2013,toselli2005domain}. Newton's method starts from whole-system linearization.
  The linearization leads to a large sparse linear system
where the matrix may be represented either explicitly by storing the nonzero coefficients or implicitly by various
``matrix-free'' approaches \cite{BrownMatrixFree1986,knollkeyes:04}. However, Newton's method has a number of drawbacks
as a stand-alone solver.  The repeated construction and solution of the linearization cause memory bandwidth and
communication bottlenecks to come to the fore with regard to performance. Also possible is a lack of convergence
robustness when the initial guess is far from the solution.  Luckily, a large design space for nonlinear solvers exists
to complement, improve, or replace Newton's method.  Only a small part of this space has yet been explored either
experimentally or theoretically.

In this paper we consider a wide collection of solvers for nonlinear equations.  In direct equivalence to the case of
linear solvers, we use a small number of algorithmic building blocks to produce a vast array of solvers with different
convergence and performance properties.  Two related solver techniques, nonlinear \composition and
preconditioning, are be used to construct these solvers. The contributions of this paper are twofold: introduction of a
systematic approach to combining nonlinear solvers mathematically and in software and demonstration of the construction
of efficient solvers for several problems of interest.  Implementations of the solvers in this paper are
available in the PETSc library and may be brought to bear on relevant applications, whether simulated on a laptop or
on a supercomputer.
\section{Background}
\label{sec:background}
We concern ourselves with the solution of nonlinear equations of the form
\begin{equation}
  \label{eq:nonprob}
  \vF(\vx) = \vb
\end{equation}
\noindent for a general discretized nonlinear function $\vF : \mathbb{R}^n \rightarrow \mathbb{R}^n$ and right hand side
(RHS) $\vb$. We define the nonlinear residual as
\begin{equation}\label{eq:residual}
  \vr(\vx) = \vF(\vx) - \vb.
\end{equation}
\noindent We retain both notations as certain solvers (\secref{sec:fas}) require the RHS to be modified. We will use 
\begin{equation}
\label{eq:jacobian}
\vJ(\vx) = \frac{\partial\vF(\vx)}{\partial\vx}
\end{equation}
to denote the Jacobian of $\vF(\vx)$.

 The linear system
\begin{equation*}
  \vA\vx = \vb
\end{equation*}
\noindent with residual
\begin{equation*}
  \vr(\vx) = \vA\vx - \vb
\end{equation*}
\noindent is an important special case from which we can derive valuable insight into solvers for the nonlinear problem.
Stationary solvers for linear systems repeatedly apply a linear operator in order to progressively improve the
solution. The application of a linear stationary solver by defect correction may be written as
\begin{equation}
 \vx_{i+1} = \vx_{i} - \vP^{-1} \left(\vA \vx_{i} - \vb\right),
\end{equation}
where $\vP^{-1}$ is a linear operator, called a preconditioner, whose action approximates, in some sense, the
inverse of $\vA$. The Jacobi, Gauss-Seidel, and multigrid iterations are all examples of linear stationary solvers.

\composition of linear preconditioners $\vP^{-1} $ and $\vQ^{-1} $ may proceed in two different ways,
producing two new stationary solvers. The first is the additive combination
\begin{eqnarray*}
\vx_{i+1} = \vx_{i} - \left(\alpha_{\vP}\vP^{-1} + \alpha_{\vQ}\vQ^{-1}\right) \left(\vA\vx_{i} - \vb\right),
\end{eqnarray*}
\noindent with weights $\alpha_{\vP}$ and $\alpha_{\vQ}$.  The second is the multiplicative combination
\begin{eqnarray*}
\vx_{i+1/2} & = & \vx_{i} - \vP^{-1}\left(\vA \vx_{i} - \vb\right) \\
\vx_{i+1} & = & \vx_{i+1/2} - \vQ^{-1}\left(\vA \vx_{i+1/2} - \vb\right). \\
\end{eqnarray*}
\compositions consisting of $\vP^{-1}$ and $\vQ^{-1}$ are an effective acceleration strategy if $\vP^{-1}$
eliminates a portion of the error space and $\vQ^{-1}$ handles the rest.  A now mature theory for these \compositions
was developed in the 1980s and 1990s in the context of domain decomposition methods \cite{smith1997domain,toselli2005domain}.

Linear left- and right-preconditioning, when used in conjunction with Krylov iterative methods
\cite{saad2003}, is standard practice for the parallel solution of linear systems of equations.  We write the use of a
linear Krylov method as $\krylov(\vA,\vx,\vb),$ where $\vA$ is the matrix, $\vx$ the initial solution, and $\vb$ the RHS.

Linear left-preconditioning recasts the problem as
\begin{equation*}
  \vP^{-1}(\vA\vx - \vb) = 0,
\end{equation*}
\noindent while right-preconditioning takes two stages,
\begin{align*}
  \vA\vP^{-1}\vy = \vb\\
  \vP^{-1}\vy = \vx,\\
\end{align*}
\noindent where one solves for the preconditioned solution $\vy$ and then transforms $\vy$ using $\vP^{-1}$
to the solution of the original problem.
\section{Nonlinear Composed Solvers}
\label{sec:composition}
We take the basic patterns from the previous section and apply them to the nonlinear case.  We emphasize
that unlike the linear case, the nonlinear case requires the approximate solution as well as the residual to be defined both in the
outer solver and in the preconditioner.  With this in mind, we will show how \composition and preconditioning may be
systematically transferred to the nonlinear case.

We use the notation $ \vx_{i+1} = \vM(\vr,\vx_{i})$ for the action of a nonlinear solver. In most cases, but not always,
$\vr = \vr(\vx_i)$, that is, $\vr $ is simply the most recent residual.  It is
also useful to consider the action of a solver that is dependent on the previous $m$ approximate solutions and the
previous $m$ residuals and as $ \vx_{i+1} = \vM(\vr,\vx_{i-m+1},\cdots,\vx_{i},\vr_{i-m+1},\cdots,\vr_{i})$.
Methods that store and use information from previous iterations such as previous solutions, residuals, or step directions will have those listed in the
per-iteration inputs as well.

Nonlinear \composition consists of a sequence or series of two (or more) solution methods $\vM$ and $\vN$, which
both provide an approximate solution to \eqnref{eq:nonprob}.  Nonlinear preconditioning, on the other hand, may be cast
as a modification of the residual $\vr$ through application of an inner method $\vN$.
The modified residual is then provided to an outer solver $\vM,$ which solves the preconditioned system.

An additive \composition may be written as
\begin{equation}
  \label{eq:additivecomposite}
  \vx_{i+1} = \vx_i + \alpha_\vM\left(\vM(\vr,\vx_{i}) - \vx_i\right) + \alpha_\vN\left(\vN(\vr,\vx_{i}) - \vx_i\right)
\end{equation}
\noindent for weights $\alpha_\vM$ and $\alpha_\vN$.  The multiplicative \composition is 
\begin{equation}
  \label{eq:multiplicativecomposite}
  \vx_{i+1} = \vM(\vr,\vN(\vr,\vx_{i})) =  \vM(\vr(\vN(\vr(\vx_{i}),\vx_{i}),\vN(\vr(\vx_{i}),\vx_{i})),
\end{equation}
which simply states: update the solution using the current solution and residual with the first solver and then update the solution again using the resulting new solution and new residual with the second solver..
\noindent Nonlinear left-preconditioning may be directly recast from the linear stationary solver case:
\begin{equation}
\vP^{-1}\left(\vA \vx - \vb\right) = 0,
\end{equation}
\noindent which we can rewrite as a fixed-point problem
\begin{equation}
  \vx - \vP^{-1}\left(\vA \vx - \vb\right) = \vx
\end{equation}
\noindent with nonlinear analog
\begin{equation}
  \vx = \vN(\vr,\vx)
\end{equation}
\noindent so the equivalent preconditioned nonlinear problem can be reached by subtracting $\vx$ from both sides, giving
\begin{equation}
  \vx - \vN(\vr,\vx) = 0.
\end{equation}
\noindent Thus, the left-preconditioned residual is given by
\begin{equation}\label{eq:leftnonlinearpreconditioning}
  \vr^l(\vx) = \vx- \vN(\vr,\vx).
\end{equation}
Under the certain circumstances the recast system will have better conditioning and less severe nonlinearities than the
original system.  The literature provides several examples of nonlinear left-preconditioning: nonlinear successive over
relaxation (\SOR) has been used in place of linear left-applied \SOR\ \cite{ChanNPNK1984}, additive
Schwarz-preconditioned inexact Newton (\ASPIN) \cite{ck02} uses an overlapping nonlinear additive Schwarz method to
provide the left-preconditioned problem for a Newton's method solver. Walker and Ni's fixed-point
preconditioned Anderson mixing \cite{walkerni} uses a similar methodology.  Many of these methods or variants thereof
are discussed and tested later in this paper.

Nonlinear right-preconditioning\footnote{Note that nonlinear right
  preconditioning is a misnomer because it does not simplify to right
  linear preconditioning in the linear case. It actually results in
  the new linear system $ A(I - P^{-1}A)y = (I - A P^{-1})b$. However,
  as with linear right preconditioning, the inner solver is applied
  before the function (matrix-vector product in the linear case)
  evaluation, hence the name.} involves a different recasting of the
residual.  This time, we treat the nonlinear solver $\vN$ as a
nonlinear transformation of the problem to one with solution $\vy$, as
$\vx = \vP^{-1}\vy$ is a linear transformation of $\vx$.  Accordingly,
the nonlinear system in \eqnref{eq:nonprob} can be rewritten as the
solution of the system perturbed by the preconditioner
\begin{equation}
  \label{eq:rightnonlinearproblem}
  \vF(\vN(\vF,\vy)) = \vb
\end{equation}
\noindent and the solution by outer solver $\vM$ as
\begin{equation}
  \label{eq:rightnonlinearpreconditioning}
  \vy_{i+1} = \vM(\vr(\vN(\vF,\cdot)),\vx_{i})
\end{equation}
\noindent followed by
\begin{equation*}
  \label{eqn:rightnonlinearfinalapplication}
  \vx_{i+1} = \vN(\vr,\vy_{i+1}).
\end{equation*}
Nonlinear right-preconditioning may be interpreted as $\vM$ putting preconditioner $\vN$ ``within striking distance'' of the solution.
Once can solve for $\vy$ in \eqnref{eq:rightnonlinearproblem} directly, with an outer solver using residual $\vr(\vN(\vr,\vx))$.  However, the
combination of inner solve and function evaluation is significantly more expensive than computing $\vF(\vx)$ and should
be avoided. We will show that when $\vM(\vr,\vx)$ is a Newton-Krylov solver,
\eqnref{eq:rightnonlinearpreconditioning} is equivalent to \eqnref{eq:multiplicativecomposite}.  Considering
them as having similar mathematical and algorithmic properties is appropriate in general.

In the special case of Newton's method, nonlinear
right-preconditioning is referred to by Cai \cite{CaiNODDM2009} and
Cai and Li \cite{CaiNKRAS2011} as nonlinear elimination
\cite{Lanzkron1996analysis}. The idea behind nonlinear elimination is
to use a local nonlinear solver to fully resolve difficult localized
nonlinearities when they begin to cause difficulties for the global
Newton's method; see \secref{sec:nk}.  Grid sequencing
\cite{KnollEnhanced1998,Smooke85GridSequence} and pseudo-transient
\cite{KelleyKeyes1998} continuation methods work by a similar
principle, using precursor solves to put Newton's method at an initial
guess from which it has fast convergence. Alternating a global linear
or nonlinear step with local nonlinear steps has been studied as the
LATIN method
\cite{Cresta2007Localization,Ladeveze2010LATIN,Ladeveze1999nonlinear}.
Full approximation scheme (\FAS) preconditioned nonlinear GMRES (\NGMRES), discussed below, for
recirculating flows \cite{ow1} is another application of nonlinear
right-preconditioning, where the \FAS\ iteration is stabilized and
accelerated by constructing a combination of several previous $\FAS$
iterates.

For solvers based on a search direction, left-preconditioning is much more natural.  In fact, general line searches
may be expressed as left-preconditioning by the nonlinear Richardson method. Left-preconditioning also has the property
that for problems with poorly scaled residuals, the inner solver may provide a tenable search direction when one could
not be found based on the original residual.
A major difficulty with the left-preconditioned option is that the function evaluation may be much more expensive.
Line searches involving the direct computation of $\vx - \vM(\vr,\vx)$ at points along the line may be overly
expensive given that the computation of the residual now requires nonlinear solves.  Line searches over $\vx -
\vM(\vr,\vx)$ may also miss stagnation of weak inner solvers and must be monitored.  One may also base the line
search on the unpreconditioned residual. The line search based on $\vx -
\vM(\vr,\vx)$ is often recommended \cite{HwangCai05} and is the ``correct'' one for general left-preconditioned
nonlinear solvers.

Our basic notation for \compositions and preconditioning is described in \tabref{tab:compprenotation}.
\begin{table}[H]
\centering
\caption{
  Nonlinear \compositions and preconditioning given outer and inner solver $\vM$ and $\vN$.
}
\label{tab:compprenotation}
\begin{tabular}{r|c|c|l}
  Composition Type              & Symbol   & Statement & Abbreviation \\
  \hline
  Additive Composite            & $+$      & $\vx + \alpha_\vM\left(\vM(\vr,\vx)-\vx\right)+\alpha_\vN\left(\vN(\vr,\vx)-\vx\right)$             & $\vM + \vN$\\
  Multiplicative Composite      & $*$      & $\vM(\vr,\vN(\vr,\vx))$                                                     & $\vM * \vN$\\
  Left-Preconditioning          & $\lp$    & $\vM(\vx - \vN(\vr,\vx),\vx)$                                               & $\vM \lp \vN$\\
  Right-Preconditioning         & $\rp$    & $\vM(\vr(\vN(\vr,\vx)),\vx)$                                                & $\vM \rp \vN$\\
  Inner Linearization Inversion & $\lin$   & $\vy = \vJ(\vx)^{-1}\vr(\vx) = \krylov(\vJ(\vx),\vy_0,\vr(\vx))$                & $\NEWT\lin\krylov$\\
\end{tabular}
\end{table}
\section{Solvers}
\label{sec:solvers}
We now introduce several algorithms, the details of their
implementation and use, and an abstract notion of how they may be
composed.  We first describe outer solution methods and how
composition is applied to them. We then move on to solvers used
primarily as inner methods.  The distinction is arbitrary but leads to
the discussion of decomposition methods in \secref{sec:decomposition}.

\subsection{Line Searches}
\label{sec:linesearch}
The most popular strategy for increasing robustness or providing globalization in the solution of nonlinear PDEs is the
line search.  Given a functional $\vf(\vx)$, a starting point $\vx_i$, and a direction $\vd$, we compute $\lambda
\approx \displaystyle\arg\min_{\mu>0}\vf(\vx_i + \mu\vd)$.  The functional $\vf(\cdot)$ may be $\|\vr(\cdot)\|_2^2$ or a
problem-specific objective function.  Theoretical guarantees of convergence may be made for many line search procedures
if $\vd$ is a descent direction.  In practice, many solvers that do not converge when only full or simply damped steps
are used converge well when combined with a line search.  Different variants of line searches are appropriate for
different solvers.  We may organize the line searches into two groups based on a single choice taken in the
algorithm: whether the full step is likely to be sufficient (for example, with Newton's method near the solution) or not.  If it is likely to be sufficient, the algorithm should default to taking the full step in a
performance-transparent way.  If not, and some shortening or lengthening of the step is assumed to be required, the
line search begins from the premise that it must determine this scaling factor.

In the numerical solution of nonlinear PDEs, we generally want to guarantee progress in the minimization of
$f(\vx)=\|\vr\|_2^2$ at each stage.  However, $\nabla\|\vr\|^2_2$ is not $\vr$, but instead $2\vJ^{\top}(\vx)\vr(\vx)$.  With
Newton's method, the Jacobian
has been computed and iteratively inverted before the line search and may be used in the Jacobian-vector product.
Also note that the transpose product is not required for the line search, since the slope of $f(\vx)$ in the direction of
step $\vy$ may be expressed as the scalar quantity $s = \vr(\vx)^T(\vJ(\vx)\vy)$.
A cubic backtracking (\BT) line search as described in Dennis and Schnabel \cite{Dennis:83} is used in conjunction with
methods based on Newton's method in this work. The one modification necessary is that it is modified to act on
the optimization problem arising from $\|\vr(\vx)\|^2_2$.  \BT\ defaults to taking the full step if that step
is sufficient with respect to the Wolfe conditions \cite{WolfeConditions}, and does no more work unless necessary.  The
backtracking line search may stagnate entirely for ill-conditioned Jacobian \cite{tuminaro2002backtracking}.
For a general step  \BT\ is not appropriate for a few reasons. First, steps arising from methods other than Newton's method are more likely to be
ill-scaled, and as such the assumption that the full-step is appropriate is generally invalid.  For this reason we only use \BT\ in the
case of Newton's method.
Second, there are also cases where we lose many of the other assumptions that make \BT\ appropriate.  Most of the
algorithms described here, for example, do not necessarily assemble the Jacobian and require many more iterations than does Newton's
method.  Jacobian assembly, just to perform the line search, in these cases would become extremely burdensome.  This situation removes any possibility of using many
of the technologies we could potentially import from optimization, including safeguards and guarantees of minimization.
However, in many cases one can still assume that the problem has optimization-like qualities.  Suppose
that $\vr(\vx)$ is the gradient of some (hidden) objective functional $\vf(\vx)$ instead of $\|\vr\|_2^2$.  Trivially,
one may minimize the hidden $\vf(\vx + \lambda\vy)$ by finding its critical points, which are roots of
\begin{equation*}
  \vy^{\top}\vr(\vx + \lambda \vy) = \frac{d\vf(\vx + \lambda \vy)}{d\lambda},
\end{equation*}
\noindent by using a secant method. The resulting critical point (\CP) line search is outlined in \agmref{alg:cplinesearch}.
\begin{algorithm}[H]
\caption{\CP\ Line Search}\label{alg:cplinesearch}
\begin{algorithmic}[1]
  \Procedure{\CP}{$\vr, \vy, \lambda_0, n$}
  \State $\lambda_{-1} = 0$
  \For{ $i = 0$ }{ $n-1$ }
  \State $\lambda_{i+1} = \lambda_{i} - \frac{\vy^{\top}\vr(\vx + \lambda_{i}\vy)(\lambda_i - \lambda_{i-1})}{\vy^{\top}\vr(\vx + \lambda_i\vy) - \vy^{\top}\vr(\vx + \lambda_{i-1}\vy)}$
  \EndFor
  \EndProcedure
  \State \Return $\lambda_n$
\end{algorithmic}
\end{algorithm}
\CP\ differs from \BT\ in that it will not minimize $\|\vr(\vx + \lambda\vy)\|_2$ over
$\lambda$.  However, \CP\ shows much more rapid convergence than does \LT\ for
certain solvers and problems.  Both line searches may be started from $\lambda_0$ as an arbitrary or problem-dependent
damping parameter, and $\lambda_{-1} = 0$. In practice, one iteration is usually satisfactory.  For highly
nonlinear problems, however, overall convergence may be accelerated by a more exact line search corresponding to a small number of
iterations.  \CP\ has also been suggested for nonlinear conjugate gradient methods~\cite{ShewchukCGWOAP1994}.

If our assumption of optimization-like qualities in the problem becomes invalid, then we can do little besides attempt
to find a minimum of $\|\vr(\vx + \lambda\vy)\|_2$. An iterative secant search for a minimum value of $\|\vr(\vx +
\lambda\vy\|_2$ is defined in \agmref{alg:l2linesearch}.

\begin{algorithm}[H]
\caption{\LT\ Line Search}\label{alg:l2linesearch}
\begin{algorithmic}[1]
        \Procedure{\LT}{$\vr, \vy, \lambda_0, n$}
        \State $\lambda_{-1} = 0$
        \For{ $i = 0$ }{ $n-1$ }
        \State $\nabla_{\vy}\|\vr(\vx + \lambda_{i}\vy)\|^2_2 = \frac{3\|\vr(\vx + \lambda_i\vy)\|^2_2
                                                                 - 4\|\vr(\vx + \frac{1}{2}(\lambda_i + \lambda_{i-1})\vy)\|_2^2
                                                                 + \|\vr(\vx + \lambda_{i-1}\vy))\|^2_2}
                                                                 {(\lambda_i - \lambda_{i-1})}$
        \State $\nabla_{\vy}\|\vr(\vx + \lambda_{i-1}\vy)\|^2_2 = \frac{\|\vr(\vx + \lambda_i\vy)\|^2_2
                                                                 - 4\|\vr(\vx + \frac{1}{2}(\lambda_i + \lambda_{i-1})\vy)\|_2^2
                                                                 + 3\|\vr(\vx + \lambda_{i-1}\vy))\|^2_2}
                                                                 {(\lambda_i - \lambda_{i-1})}$
        \State $\lambda_{i+1} = \lambda_{i} - \frac{\nabla_{\vy}\|\vr(\vx + \lambda_{i}\vy)\|^2_2(\lambda_i - \lambda_{i-1})}
               {\nabla_{\vy}\|\vr(\vx + \lambda_i\vy)\|_2^2 - \nabla_{\vy}\|\vr(\vx + \lambda_{i-1}\vy)\|_2^2}$
        \EndFor
        \EndProcedure
        \State\Return $\lambda_n$
\end{algorithmic}
\end{algorithm}
When converged, \LT\ is equivalent to an optimal damping in the direction of the residual and will forestall
divergence.  $\nabla_{\vy}$ is calculated by polynomial approximation, requiring two additional residual evaluations per
application.   In practice and in our numerical experiments the number of inner iterations $n$ is $1$.

\subsection{Nonlinear Richardson (\NRICH)}
\label{sec:nrich}
The nonlinear analogue to the Richardson iteration is merely the simple application of a line search.  \NRICH\ takes a
step in the negative residual direction and scales that step sufficiently to guarantee convergence.  \NRICH\ is known as
steepest descent~\cite{goldsteinsteepest} in the optimization context where $\vr$ is the gradient of a functional
$\vf(\cdot)$ noted above.  \NRICH\ is outlined in \agmref{alg:nrich}.
\begin{algorithm}[H]
\caption{Nonlinear Richardson Iteration}\label{alg:nrich}
\begin{algorithmic}[1]
        \Procedure{\NRICH}{$\vr,\vx_i$}
        \State $\phantom{\vx_{i+1}}
        \mathllap{\vd} = -\vr(\vx_i)$
        \State $\vx_{i+1} = \vx_i + \lambda\vd$ \Comment{$\lambda$ determined by line search}
        \EndProcedure
        \State \Return $\vx_{i+1}$
\end{algorithmic}
\end{algorithm}

\NRICH\ is often slow to converge for general problems and stagnates quickly.  However, different step
directions than $\vr$ may be generated by nonlinear preconditioning and can improve convergence dramatically.
\begin{algorithm}[H]
\caption{Nonlinear Richardson Iteration: Left-Preconditioned by $\vM()$}\label{alg:nrich:left}
\begin{algorithmic}[1]
        \Procedure{\NRICH}{$\vx - \vM(\vr,\vx), \vx_i$}
        \State $\phantom{\vx_{i+1}}\mathllap{\vd} = \vM(\vr,\vx_i) - \vx_i$
        \State $\vx_{i+1} = \vx_i + \lambda\vd$ \Comment{$\lambda$ determined by line search}
        \EndProcedure
        \State\Return $\vx_{i+1}$
\end{algorithmic}
\end{algorithm}
As shown in \agmref{alg:nrich:left}, we replace the original residual equation $\vr(\vx)$ with $\vx-\vM(\vr,\vx)$
and apply \NRICH\ to the new problem. There are two choices for $\vr$ in the line search. The
first is based on minimizing the original residual, as in the unpreconditioned case; the second minimizes the norm
of the preconditioned residual instead.  Minimizing the unpreconditioned residual with a preconditioned step is more
likely to stagnate, as there is no guarantee that the preconditioned step is a descent direction
with respect to the gradient of $\|\vr(\vx)\|_2^2$.
\subsection{Anderson Mixing (\ANDERSON)}
\label{sec:anderson}
\ANDERSON  \cite{AndersonMixing} constructs a new approximate solution as a
combination of several previous approximate solutions and a new trial.

\begin{algorithm}[H]
\caption{Anderson Mixing}\label{alg:anderson}
\begin{algorithmic}[1]
        \Procedure{$\ANDERSON$}{$\vr,\vx_i \cdots \vx_{i-m+1}$ }
        \State $\vx_i^M = \vx_i + \lambda\vr(\vx_i)$
        \State \textbf{minimize} $\left\|\vr(\left(1 - \sum^{i-1}_{k=i-m}\alpha_k\right)\vx_i^M + \sum_{k = i - m}^{i-1}\alpha_k\vx_k)\right\|_2$ over $\left\{\alpha_{i - m} \cdots \alpha_{i}\right\}$
        \State $\vx_{i+1} = \left(1 - \sum^{i}_{k=i-m}\alpha_k\right)\vx_i^M + \sum^{i-1}_{k = i - m}\alpha_k\vx^M_k$
        \EndProcedure
        \State \Return $\vx_{i+1}$
\end{algorithmic}
\end{algorithm}
In practice, the nonlinear minimization problem in \agmref{alg:anderson} is simplified.  The $\alpha_i$ are
computed by considering the linearization
\begin{equation}
  \label{eq:ngmreslinearization}
  \vr((1 - \sum^{i-1}_{k=i-m}\alpha_k)\vx^M_i + \sum^{i-1}_{k=i-m}\vx^M_k) \approx (1 - \sum^{i-1}_{k=i-m}\alpha_k )\vr(\vx^M_i) + \sum^{i-1}_{k=i-m}\vr(\vx^M_k)
\end{equation}
\noindent and solving the related linear least squares problem
\begin{equation}
  \label{eq:minprob}
  [\vr(\vx_j) - \vr(\vx^M)]^{\top}[\vr(\vx_i) - \vr(\vx^M)]\alpha_j = \vr(\vx_i)^{\top}[\vr(\vx_i) - \vr(\vx^M)].
\end{equation}
\ANDERSON\ solves \eqnref{eq:minprob} by dense factorization. The work presented here uses the SVD-based least-squares
solve from LAPACK in order to allow for potential singularity of the system arising from stagnation to be detected
outside of the subsystem solve.
A number of related methods fall under the broad category of nonlinear series accelerator methods.  These include
Anderson mixing as stated above. Nonlinear GMRES \NGMRES\ \cite{washio1997krylov} is a variant that includes conditions to avoid stagnation
and direct inversion in the iterative subspace (DIIS) \cite{PulayCA1980}, which formulates the minimization problem in an alternative fashion.
Both right- and left-preconditioning can be applied to Anderson mixing.  With right-preconditioning the computation of
$\vx_i^M = \vx_{i-1} + \lambda\vd$ is replaced with $\vx^M = \vM(\vr,\vx_{i-1})$ and $\vr(\vx_i)$ with $\vr(\vx^M)$.
This incurs an extra function evaluation, although this is usually included in the nonlinear preconditioner application.
Left-preconditioning can be applied by replacing $\vr$ with $\vr^l$.
Right-preconditioned \NGMRES\ has been applied for recirculating flows \cite{ow1} using \FAS\, and in stabilizing lagged
Newton's methods \cite{CarlsonGWMFE1998, Scott2003Krylov}.  Simple preconditioning of \NGMRES\ has also been proposed in
the optimization context \cite{sterck2012steepest}.  Preconditioned Anderson mixing has been leveraged as an outer
accelerator for the Picard iteration \cite{Lott2012Accelerated,walker2010accelerated}.

We can use the same formulation expressed in \eqnref{eq:minprob} to determine $\alpha_\vM$ and $\alpha_\vN$ (or any
number of weights) in \eqnref{eq:additivecomposite}, using the solutions and final residuals from a series of
inner nonlinear solvers instead of the sequence of previous solutions.  This sort of residual-minimizing technique
is generally applicable in additive \compositions of solvers.  All instances of additive
\composition in the experiments use this formulation.
\subsection{Newton-Krylov Methods (\NK)}
\label{sec:nk}
In Newton-Krylov methods, the search direction is determined by inexact iterative inversion of the Jacobian
applied to the residual by using a preconditioned Krylov method.
\begin{algorithm}[H]
\caption{Newton-Krylov Method}\label{alg:nk}
\begin{algorithmic}[1]
  \Procedure{\NK}{$\vr,\vx_i$}
  \State $\phantom{\vx_{i+1}}\mathllap{\vd} = \vJ(\vx_i)^{-1} \vr(\vx_i)$ \Comment{approximate inversion by Krylov method}
  \State $\vx_{i+1} = \vx_i + \lambda \vd$ \Comment{$\lambda$ determined by line search}
  \EndProcedure
  \State \Return $\vx_{i+1}$
\end{algorithmic}
\end{algorithm}
\NK\ cannot be guaranteed to converge far away from the solution, and it is routinely enhanced by a line search.
 In our formalism, \NK\ with a line search can be expressed as \NRICH\ left-preconditioned
by \NK\, coming from \agmref{alg:nrich:left}.  However, we will consider line-search globalized \NK\ as the standard,
and will omit the outer \NRICH\ when using it, leading to \agmref{alg:nk}.  Other globalizations, such as trust-region
methods, are also often used in the context of optimization but will not be covered here.

\begin{figure}[H]
\centering
    \includegraphics[width=0.8\textwidth]{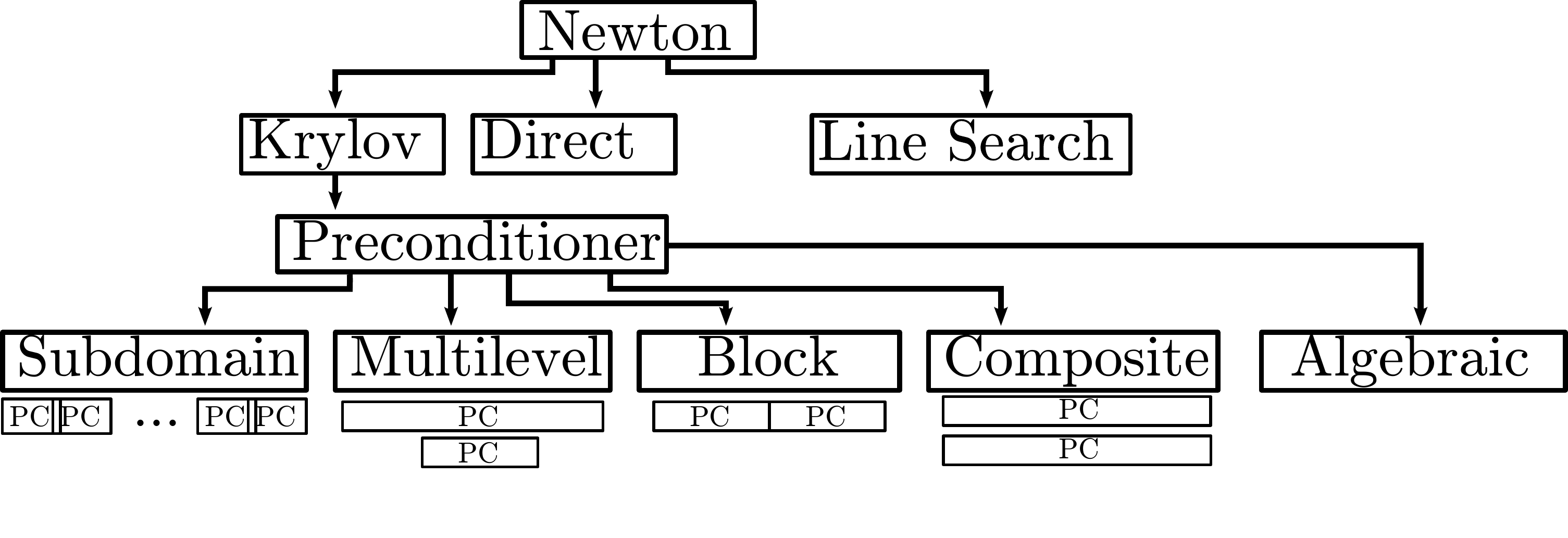}
    \caption{Organization of a Newton-Krylov solver with potential points of customization.
      Although one can alter the type, parameters, or tolerances of Newton's
      method, the Krylov solver, or the line search, we see that the majority of the potential customization and composition
      occurs at the level of the preconditioner.  Possibilities include the use of algebraic preconditioners, such as direct solvers
      or sparse factorizations, or domain decomposition, multigrid, or composite solvers, which all have subcomponent
      solvers and preconditioners. These divisions may either be by subdomain, as with the additive Schwarz methods, or by physics field components, as with block preconditioners (for example, one block for pressures and one for velocities).}
    \label{fig:nkdiagram}
\end{figure}

\NK\ is the general-purpose workhorse of a majority of simulations requiring solution of nonlinear equations.
Numerous implementations and variants exist, both in the literature and as software \cite{Brown1994Convergence,
  Eisenstat1994Globally}.  The general organization of the components of \NK\ is shown in \figref{fig:nkdiagram}.  Note
that the vast majority of potential customization occurs at the level of the linear preconditioner and that access to
fast solvers is limited to that step.  
Nonlinear right-preconditioning of Newton's method may take two forms.  We describe our choice and comment on the
alternative.  At first glance, right-preconditioning of Newton's method requires construction or application
of the right-preconditioned Jacobian
\begin{equation}
\label{eq:rightpcjac}
\frac{\partial\vF(\vM(\vr,\vy))}{\partial\vy} = \vJ(\vM(\vr,\vy))\frac{\partial\vM(\vr,\vy)}{\partial\vy},
\end{equation}
\noindent and we note that the computation of $\frac{\partial\vM(\vr,\vy_i)}{\partial \vy_i}$ is
generally impractical.  Consider the application of nonlinear right preconditioning on an iteration of Newton's method, which would proceed in two steps:
\begin{align*}
  \vy_{i+1}               &= \vx_i - \lambda\left(\frac{\partial\vM(\vr,\vx_i)}{\partial\vx_i}\right)^{-1}
                                    \vJ\left(\vM(\vr,\vx_i)\right)^{-1}\vr(\vM(\vr,\vx_i)) \\
  \vx_{i+1}               &= \vM(\vr,\vy_{i+1}).
\end{align*}
We may nest these two steps as
\begin{equation}
  \vx_{i+1}  = \vM(\vr,\vx_i - \lambda\left(\frac{\partial\vM(\vr,\vx_i)}{\partial\vx_i}\right)^{-1}\vJ(\vM(\vr,\vx_i))^{-1}\vr(\vM(\vr,\vx_i))),
\end{equation}
and by Taylor expansion one may say that
\begin{equation*}
\vM(\vr,\vx - \lambda\vd) = \vM(\vr,\vx) - \lambda\left(\frac{\partial\vM(\vr,\vx)}{\partial\vx}\right)\vd + ...,
\end{equation*}
\noindent giving us the simplification
\begin{align*}
x_{i+1}      &  \approx \vM(\vr,\vx_i) -
                            \lambda\left(\frac{\partial\vM(\vr,\vx_i)}{\partial\vx_i}\right)
                            \left(\frac{\partial\vM(\vr,\vx_i)}{\partial\vx_i}\right)^{-1}
                            \vJ(\vM(\vr,\vx_i))^{-1}\vr(\vM(\vr,\vx_i)))\\
                         &= \vM(\vr,\vx_i) - \lambda\vJ(\vM(\vr,\vx_i))^{-1}\vr(\vM(\vr,\vx_i))\\
\end{align*}
In this form, \NK$(\vF(\vM(\vr,\vx)),\vx)$ is easily implemented, as shown in \agmref{alg:nkmright}.
\begin{algorithm}[H]
\caption{Right-Preconditioned Newton-Krylov Method}
\label{alg:nkmright}
\begin{algorithmic}[1]
        \Procedure{NK}{$\vr(\vM(\vr,\vx)),\vx_i$}
        \State $\phantom{\vx_{i+1}}
        \mathllap{\vx_{i+\frac{1}{2}}} = \vM(\vr,\vx_i)$
        \State $\phantom{\vx_{i+1}}
        \mathllap{\vd} = \vJ(\vx_{i+\frac{1}{2}})^{-1}\vr(\vx_{i+\frac{1}{2}})$
        \State $\vx_{i+1} = \vx_{i+\frac{1}{2}} + \lambda \vd$  \Comment{$\lambda$ determined by line search}
        \EndProcedure
        \State \Return $\vx_{i+1}$
\end{algorithmic}
\end{algorithm}
Note that $\NK \rp \vM$ is merely $\vM * \NK$ in this form: one solver runs after the other.  The connection between
composition and preconditioning in Newton's method provides an inkling of the advantages afforded to \NK\ as an inner
composite solver combined with other methods. An alternative approach applies an approximation to \eqnref{eq:rightpcjac}
directly~\cite{Birken2009NonPreNK}.  The approximation is
\begin{align*}
  \vF(\vM(\vr,\vy_i)) &=
                             \vJ(\vM(\vr,\vy_i))
                             \left(\frac{\partial\vM(\vr,\vy_i)}{\partial\vy_i}\right)
                             (\vy_{i+1} - \vy_i)\\
                          &\approx
                             \vJ(\vM(\vr,\vy_i))
                             (\vM(\vr,\vy_i + [\vy_{i+1} - \vy_{i}]) - \vx_i),
\end{align*}
\noindent which is solved for $\left[\vy_{i+1} - \vy_{i}\right]$.  Right-preconditioning by the approximation requires an application of the nonlinear
preconditioner for every inner Krylov iterate and limits our choices of Krylov solver to those tolerant of
nonlinearity, such as flexible GMRES \cite{Saad1993}.

\begin{algorithm}[H]
\caption{Left-Preconditioned Newton-Krylov Method}\label{alg:nkleft}
\begin{algorithmic}[1]
        \Procedure{\NK}{$\vx - \vM(\vr,\vx),\vx_i$}
        \State $\phantom{\vx_{i+1}}\mathllap{\vd} = \frac{\partial(\vx_i - \vM(\vr,\vx_i))}{\partial\vx_i}^{-1} (\vx_i - \vM(\vr,\vx_i))$ \Comment{approximate inversion by Krylov method}
        \State $\vx_{i+1} = \vx_i + \lambda \vd$  \Comment{$\lambda$ determined by line search}
        \EndProcedure
        \State \Return $\vx_{i+1}$
\end{algorithmic}
\end{algorithm}

For nonlinear left-preconditioning, shown in \agmref{alg:nkleft}, we replace the computation $\vr(\vx_i)$ with $\vx_i -
\vM(\vr,\vx_i)$ and take the Jacobian of the function as an approximation of
\begin{equation}
  \frac{\partial(\vx_i - \vM(\vr,\vx_i))}{\vx_i} = \vI - \frac{\partial\vM(\vr,\vx_i)}{\partial\vx_i},
\end{equation}
\noindent where now the Jacobian is a linearization of $\vM(\vr,\vx)$ and impractical to compute in most cases.  In
the particular case where the preconditioner is the nonlinear additive Schwarz method (\NASM) (\secref{sec:nasm}), this is known as \ASPIN.  In the case of
\NASM\ preconditioning, one has local block Jacobians, so the approximation of the preconditioned Jacobian as
\begin{equation}
  \frac{\partial(\vx - \vM(\vr,\vx))}{\partial\vx} = \frac{\partial(\vx - (\vx - \sum_b\vJ^b(\vx^b)^{-1}
    \vr^b(\vx^b)))}{\partial\vx} \approx \sum_b\vJ^b(\vx^{b*})^{-1}\vJ(\vx)
\end{equation}
\noindent is used instead.  The iteration requires only one inner nonlinear iteration and some small number of block
solves, per outer nonlinear iteration. By contrast, direct differencing would require one inner iteration at each inner
linear iteration for the purpose of repeated approximate Jacobian application.  Note the similarity between \ASPIN\ and
the alternative form of right-preconditioning in terms of required components, with \ASPIN\ being more
convenient in the end.
\subsection{Quasi-Newton (\QN)}
\label{sec:qn}
The class of methods for nonlinear systems of equations known as quasi-Newton methods \cite{DennisQN} uses
previous changes in residual and solution or other low-rank expansions \cite{klie2006nonlinear} to form an approximation
of the inverse Jacobian by a series of low rank updates. The approximate Jacobian inverse is then used to compute the
search direction. If the update is done directly, it requires storing and updating the dense matrix $\vK \approx
\vJ^{-1}$, which is impractical for large systems.  One may overcome this limitation by using limited-memory variants of
the method \cite{liu1989limited,NocedalLQN,ByrdQNLimited1994,MatthiesNonlinearFEM1979}, which apply the approximate inverse Jacobian by a series
of low-rank updates.  The general form of \QN\ methods is similar to that of \NK\ and is shown in \agmref{alg:qn}.
\begin{algorithm}[H]
\caption{Quasi-Newton Update}\label{alg:qn}
\begin{algorithmic}[1]
  \Procedure{$\QN$}{$\vr, \vx_i, \vs_{i-1} \cdots \vs_{i-m}, \vy_{i-1} \cdots \vy_{i-m}$}
  \State $\vx_{i+1} = \vx_{i} - \lambda\vK_i(\vK_0,\vs_{i-1} ... \vs_{i-m}, \vy_{i-1} .. \vy_{i-m})\vr(\vx_i)$ \Comment{$\lambda$ determined by line search}
  \State $\phantom{\vx_{i+1}}
  \mathllap{\vs_{i}} = \vx_{i+1} - \vx_i$
  \State $\phantom{\vx_{i+1}} \mathllap{\vy_{i}} = \vr(\vx_{i+1}) - \vr(\vx_i)$
  \EndProcedure
  \State \Return $\vx_{i+1}$
\end{algorithmic}
\end{algorithm}
A popular variant of the method, L-BFGS, can be applied efficiently by a two-loop recursion \cite{NocedalLQN}:
\begin{algorithm}[H]
\caption{Two-Loop Recursion for L-BFGS}\label{alg:2lr}
\begin{algorithmic}[1]
  \Procedure{$\vK_i$}{$\vK_{0},\vs_{i-1} ... \vs_{i-m},\vy_{i-1} ... \vy_{i-m},\vr(\vx_i)$}
  \State $\vd_1 = \vr(\vx_i)$
  \For{ $k=i-1$ }{ $i-m$ }
  \State $\alpha_k = \frac{\vs_k^\top\vd_1}{\vy_k^\top\vs_k}$
  \State $\vd_1 = \vd_1 - \alpha_k\vy_k$
  \EndFor
  \State $\vd_2 = \vK_{0}\vd_1$
  \For{ $k=i-m$ }{ $i-1$ }
  \State $\beta_k = \frac{\vy_k^\top\vd_2}{\vy_k^\top\vs_k}$
  \State $\vd_2 = \vd_2 + (\alpha_k - \beta_k)\vs_k$
  \EndFor
  \EndProcedure
  \State \Return $\vd_2$
\end{algorithmic}
\end{algorithm}

Note that the update to the Jacobian in \agmref{alg:2lr} is symmetric and shares many of the same limitations as nonlinear conjugate gradients (\NCG)
(\secref{sec:ncg}).  It may be perplexing, from an optimization perspective, to see L-BFGS and \NCG\ applied to
nonlinear PDEs as opposed to optimization problems.  However, their exclusion from our discussion
would be as glaring as leaving conjugate gradient methods out of a discussion of Krylov solvers for linear PDEs.  In
this paper we adhere to the limitations of \NCG\ and \QN\ and use them only for problems with symmetric Jacobians.  In the case where L-BFGS is inapplicable, Broyden's ``good'' and
``bad'' methods \cite{broyden1965class} have efficient limited memory constructions
\cite{ByrdQNLimited1994,deuflhard1990fast,fangsaad} and are recommended instead \cite{NW99}. The overall performance of
all these methods, however, may rival that of Newton's method \cite{Kelley95} in many cases.  \QN\ methods also have
deep connections with Anderson mixing \cite{eyert}.  Anderson mixing, and the Broyden methods belong to a general class
of Broyden-like methods \cite{fangsaad}.

The initial approximate inverse Jacobian $\vK_{0}$ is often taken as the weighted \cite{ShannoScaling1970} identity.
However, one also can take $\vK_{0} = \vJ_{0}^{-1}$ for lagged Jacobian $\vJ_{0}$.  Such schemes greatly
increase the robustness of lagged Newton's methods when solving nonlinear PDEs \cite{brown2013quasinewton}.
Left-preconditioning for \QN\ is trivial, and either the preconditioned or unpreconditioned residual may be used with
the line search.
\subsection{Nonlinear Conjugate Gradients}
\label{sec:ncg}
Nonlinear conjugate gradient methods \cite{FletcherReeves1964} are a simple extension of the linear \CG\ method but
with the optimal step length in the conjugate direction determined by a line search rather than the exact formula that
works in the linear case.  \NCG\ is outlined in \agmref{alg:ncg}.  \NCG\ requires storage of one additional vector for
the conjugate direction $\vc_i$ but has significantly faster convergence than does \NRICH\ for many problems.
\begin{algorithm}[H]
\caption{Nonlinear CG}\label{alg:ncg}
\begin{algorithmic}[1]
        \Procedure{\NCG}{$\vr, \vx_{i},\vc_{i-1},\vr_{i-1}$}
        \State $\phantom{\vx_{i+1}}
        \mathllap{\vr_i} = \vr(\vx_i)$
        \State $\phantom{\vx_{i+1}}
        \mathllap{\beta_{i}} = \frac{\vr_{i}^{\top}(\vr_{i} - \vr_{i-1})}{\vr_{i-1}^\top\vr_{i-1}}$
        \State $\phantom{\vx_{i+1}}
        \mathllap{\vc_{i}} = -\vr(\vx_{i}) + \beta_{i}\vc_{i-1}$
        \State $\vx_{i+1} = \vx_{i} + \lambda\vc_{i}$ \Comment{$\lambda$ determined by line search}
        \EndProcedure
        \State \Return $\vx_{i+1}$
\end{algorithmic}
\end{algorithm}
Several choices exist for constructing the parameter $\beta_{i}$
\cite{daiyuan,Fletcher1987,FletcherReeves1964,hs:52}.  We choose the Polak-Ribi\`ere-Polyak \cite{Polak1969} variant.
The application of nonlinear left-preconditioning is straightforward.  Right-preconditioning is not conveniently
applicable as is the case with \QN\ in \secref{sec:qn} and linear CG.

\NCG\ has limited applicability because it suffers from the same issues as its linear cousin for problems with non-symmetric
Jacobian.  A common practice when solving nonlinear PDEs with \NCG\ is to rephrase the problem in terms of
the normal equations, which involves finding the root of $\vr_N(\vx)= \vJ^\top \vr(\vx)$.  Since the Jacobian must be
assembled and multiplied at every iteration, however, the normal equation solver is not a reasonable algorithmic choice, even
with drastic lagging of Jacobian assembly. The normal equations also have a much worse condition number than does the
original PDE. In our experiments, we use the conjugacy-ensuring \CP\ line search.
\section{Decomposition Solvers}
\label{sec:decomposition}
An extremely important class of methods for linear problems is based on domain or hierarchical decomposition, and so 
we consider nonlinear variants of domain decomposition and multilevel algorithms.  We introduce three different nonlinear solver algorithms
based on point-block solves, local subdomain solves, and coarse-grid solves.

These methods require a trade-off between the amount of computation dedicated to solving local or coarse problems and the
communication for global assembly and convergence monitoring.  Undersolving the subproblems wastes the
communication overhead of the outer iteration, and oversolving exacerbates issues of load imbalance and quickly has
diminishing returns with respect to convergence.  The solvers in this section exhibit a variety of features related to these
trade-offs.

The decomposition solvers do not guarantee convergence.  While they might converge,
the obvious extension is to use them in conjunction with the global solvers as nonlinear preconditioners or
accelerators.  As the decomposition solvers expose more possibilities for parallelism or acceleration, their effective use in
the nonlinear context should provide similar benefits to the analogous solvers in the linear context.  A major
disadvantage of these solvers is that each of them requires additional information about the local problems, a
decomposition, or hierarchy of discretizations of the domain.
\subsection{Gauss-Seidel-Newton (\GSN)}
\label{sec:gsn}
Effective methods may be constructed by exact solves on subproblems.  Suppose that the problem easily
decomposes into $n_b$ small block subproblems. If Newton's method is applied multiplicatively by blocks, the resulting
algorithm is known as Gauss-Seidel-Newton and is shown in \agmref{alg:gsn}.  Similar methods are commonly used as
nonlinear smoothers \cite{HackbuschNonlinearMG,HensonNMGOverview}. To construct \GSN, we define the individual block
Jacobians $\vJ^b(\vx^b)$ and residuals $\vr^b(\vx^b)$; $\rsrct^b$, which restricts to a block; and $\intrp^b$, which
injects the solution to that block back into the overall solution.  The point-block solver runs until the norm of the
block residual is less than $\epsilon^b$ or $m_b$ steps have been taken.
\begin{algorithm}[H]
\caption{Gauss-Seidel-Newton}\label{alg:gsn}
\begin{algorithmic}[1]
  \Procedure{\GSN}{$\vr, \vx_i$}
        \State $\vx_{i+1} = \vx_{i}$
        \For{ $b = 1$ }{ $n_{b}$ }
        \State $\vx^b_{i,0} = \rsrct^b\vx_{i+1}$
        \While{ $\|\vr^b\|_2 > \epsilon^b$ and $j < m_b$ }
        \State $\phantom{\vx_{i,j}^b}
        \mathllap{j} = j+1$
        \State $\vx_{i,j}^b = \vx_{i,j-1}^b - \vJ^b(\vx_{i,j-1}^b)^{-1}\vr_{i,j-1}^b$  \Comment{direct inversion of block Jacobian}
        \State $\vr_{i,j}^b = \vr^b(\vx^b_{i,j})$
        \EndWhile
        \State $\vx_{i+1} = \vx_{i+1} - \intrp^b(\vx_{i,0}^b - \vx_{i,j}^b)$
        \EndFor
        \EndProcedure
        \State \Return $\vx_{i+1}$
\end{algorithmic}
\end{algorithm}
\GSN\ solves small subproblems with a fair amount of computational work per degree of freedom and thus has high arithmetic
intensity.  It is typically greater than two times the work per degree of freedom of \NRICH\ for scalar problems, and even more for
vector problems where the local Newton solve is over several local degrees of freedom.

When used as a solver, \GSN\ requires one reduction per iteration.  However, stationary solvers, in both the linear and
nonlinear case, do not converge robustly for large problems.  Accordingly, \GSN\ would typically be used as a
preconditioner, making monitoring of global convergence, within the GSN iterations, unnecessary. Thus, in practice there is no synchronization;
\GSN\ is most easily implemented with multiplicative update on serial subproblems and additively in parallel.
\subsection{Nonlinear Additive Schwarz Method}
\label{sec:nasm}
\GSN\ can be an efficient underlying kernel for sequential nonlinear solvers. For parallel computing, however, an
additive method with overlapping subproblems has desirable properties with respect to communication and arithmetic
intensity.  The nonlinear additive Schwarz methods allow for medium-sized subproblems to be solved with a general
method, and the corrections from that method summed into a global search direction.  Here we limit ourselves to
decomposition by subdomain rather than splitting the problem into fields.  While eliminating fields might be tempting,
the construction of effective methods of this sort is significantly more involved than in the linear case, and many
interesting problems do not have such a decomposition readily available.  Using subdomain problems $\vF^B$, restrictions
$\rsrct^B$, injections $\intrp^B$, and solvers $\vM^B$, \agmref{alg:nasm} outlines the \NASM\ solver with $n_B$ subdomains.
\begin{algorithm}[H]
\caption{Nonlinear Additive Schwarz}\label{alg:nasm}
\begin{algorithmic}[1]
        \Procedure{\NASM}{$\vr, \vx_i$}
        \For{ $B = 1$ }{ $n_{B}$ }
        \State $\vx^B_{0}= \rsrct^B\vx_i$
        \State $\phantom{\vx^B_{0}}
        \mathllap{\vx^B} = \vM^B(\vr^B,\vx^B_{0})$
        \State $\phantom{\vx^B_{0}}
        \mathllap{\vy^B} = \vx^B_{0} - \vx^B$
        \EndFor
        \State $\vx_{i+1} = \vx_{i} - \displaystyle\sum_{B=1}^{n_B}\intrp^B\vy^B$
        \EndProcedure
        \State \Return $\vx_{i+1}$
\end{algorithmic}
\end{algorithm}
Two choices exist for the injections $\vP^B$ in the overlapping regions. We will use \NASM\ to denote the variant
that uses overlapping injection corresponding to the whole subproblem step.  The second variant, restricted additive Schwarz
(\RAS), injects the step in a non-overlapping fashion.  The subdomain solvers are typically Newton-Krylov, but the other
methods described in this paper are also applicable.
\subsection{Full Approximation Scheme}
\label{sec:fas}
\FAS\ \cite{brandt1977multi} accelerates convergence by advancing the nonlinear solution on a series of coarse rediscretizations of the
problem.  As with standard linear multigrid, the cycle may be constructed either additively or multiplicatively, with
multiplicative being more effective in terms of per-iteration convergence.  Additive approaches provide the usual advantage of
allowing the coarse-grid corrections to be computed in parallel.  We do not consider additive \FAS\ in this paper.

Given the smoother $\vM_s(\vr, \vx)$ at each level, as well as restriction ($\rsrct$), prolongation ($\intrp$) and
injection ($\inj$) operators, and the coarse nonlinear function $\vF^H$, the \FAS\ V-cycle takes the form shown in
\agmref{alg:fas}.
\begin{algorithm}[H]
\caption{Full Approximation Scheme}\label{alg:fas}
\begin{algorithmic}[1]
        \Procedure{$\FAS$}{$\vr, \vx_i$}
        \State $\phantom{\vx_{i+1}}
        \mathllap{\vx_s} = \vM_s(\vr, \vx_i)$
        \State $\phantom{\vx_{i+1}}
        \mathllap{\vx^H_i} = \inj \vx_s$
        \State $\phantom{\vx_{i+1}}
        \mathllap{\vb^H} = \rsrct[\vb  - \vF(\vx_s)] + \vF^H(\vx^H_i)$
        \State $\phantom{\vx_{i+1}}
        \mathllap{\vx_c} = \vx_s + \intrp[\FAS(\vF^H - \vb^H,\vx^H_i) - \vx^H_i]$
        \State $\phantom{\vx_{i+1}}
        \mathllap{\vx_{i+1}} = \vM_s(\vr, \vx_c)$
        \EndProcedure
        \State \Return $\vx_{i+1}$
\end{algorithmic}
\end{algorithm}
The difference between \FAS\ and linear multigrid is the construction of the coarse RHS $\vb^H$.  In \FAS\, it is guaranteed
that if an exact solution $\vx^*$ to the fine problem is found,
\begin{equation}
  \FAS(\vr^H - \vb^H,\inj\vx^*) - \inj\vx^* = 0.
\end{equation}
\noindent The correction is not necessary in the case of linear multigrid.  However, \FAS\ is mathematically identical
to standard multigrid when applied to a linear problem.

The stellar algorithmic performance of linear multigrid methods is well documented, but the arithmetic intensity may be
low. \FAS\ presents an interesting alternative to $\NK\pc\MG$, since it may be configured with high arithmetic
intensity operations at all levels. The smoothers are themselves nonlinear solution methods. \FAS-type
methods are applicable to optimization problems as well as nonlinear PDEs \cite{brandt2003multigrid}, with optimization
methods used as smoothers \cite{nash2000multigrid}.
\subsection{Summary}
\begin{figure}[H]
\centering
\includegraphics[width=\textwidth]{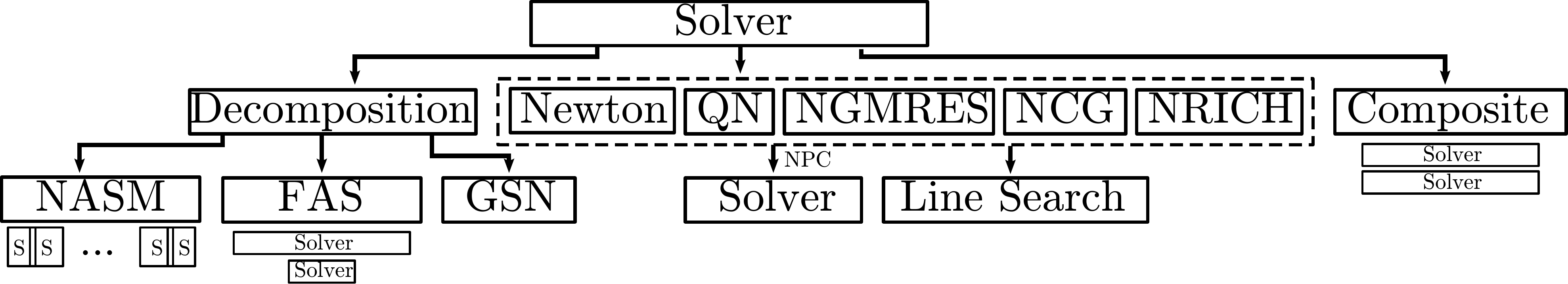}
\caption{Organization of the components of a composed nonlinear solver.  We discard the difference between solvers and
  preconditioners and see that the nesting and potential for recursive customization lives at every level of the tree.
  Possibilities include iterative solvers (including Newton's method) with nonlinear preconditioning, composite solvers
  consisting of several subsolvers, or decomposition solvers consisting of subdomain or coarse nonlinear solvers.}
\label{fig:npcdiagram}
\end{figure}

We have now introduced the mathematical construction of nonlinear composed solvers.  We have also described two classes
of nonlinear solvers, based on solving either the global problem or some partition of it.  The next step is to
describe a set of test problems with different limitations with respect to which methods will work for them, construct a
series of instructive example solvers using composition of the previously defined solvers, and test them. We have made
an effort to create flexible and robust software for composed nonlinear solvers.  The organization of the software is in
the spirit of PETSc and includes several interchangeable component solvers, including \NK; iterative solvers such as
\ANDERSON\ and \QN\ that may contain an inner preconditioner; decomposition solvers such as \NASM\ and \FAS, which are built
out of subdomain solvers; and metasolvers implementing \compositions.  A diagram enumerating the composed
solver framework analogous to that of the preconditioned \NK\ case is shown in \figref{fig:npcdiagram}.
\section{Experiments}
\label{sec:experiments}

We will demonstrate the efficacy nonlinear composition and preconditioning by experiments with a suite of nonlinear
partial differential equations that show interesting behavior in regimes with difficult nonlinearities.  These problems
are nonlinear elasticity, the lid-driven cavity with buoyancy, and the $p$-Laplacian.

One goal of this paper is to be instructive.  We first try to solve the problem efficiently with the standard solvers,
then choose a subset of the above solvers for each problem and show how to gain advantage by using composition
methods.  Limitations of discretization or problem regime are noted, and we discuss how the solvers may be used under
these limitations.  We also explain how readers may run the examples in this paper and experiment with the solvers
both on the test examples shown here and on their own problems. We feel we have no ``skin in the game'' and are not trying to show that some approaches are better or worse than others.
\subsection{Methods}
\label{sec:methods}
Our set of test algorithms is defined by the preconditioning of an iterative solver with a decomposition solver,
or the \composition of two solvers with different advantages.  In the case of the decomposition solvers, we choose one or two configurations
per test problem, in order to avoid the combinatorial explosion of potential methods that already is apparent from the two
forms of composition combined with the multitude of methods.

Our primary measure of performance is time to solution, which is both problem and equipment dependent.  Since it depends
strongly on the relative cost of function evaluation, Jacobian assembly, matrix multiplication, and subproblem or coarse
solve, we also record number of nonlinear iterations, linear iterations, linear preconditioner applications, function
evaluations, Jacobian evaluations, and nonlinear preconditioner applications.  These measures allow us to characterize
efficiency disparities in terms of major units of computational work.

We make a good faith effort to tune each method for performance while keeping subsolver parameters invariant.
Occasionally, we err on the side of oversolving for the inner solvers, to make the apples-to-apples comparison more
consistent between closely related methods.  The results here should be taken as a rough guide to improving solver
performance and robustness.  The test machine is a 64-core AMD Opteron 6274-based~\cite{amd6274} machine with 128 GB of
memory.  The problems are run with one MPI process per core.  Plots of convergence are generated with Matplotlib
\cite{Hunter2007matplotlib}, and plots of solutions are generated with Mayavi \cite{ramachandran2011mayavi}.
\subsection{Nonlinear Elasticity}
\label{sec:elasticity}
A Galerkin formulation for nonlinear elasticity may be stated as
\begin{equation}
  \label{eq:elasticity}
  \int_{\Omega} F\cdot S : \nabla v\,d\Omega + \int_{\Omega}b\cdot v\,d\Omega = 0
\end{equation}
\noindent for all test functions $v \in \mathcal{V}$; $F = \nabla \vu + I$ is the deformation gradient.  We use the Saint
Venant-Kirchhoff model of hyperelasticity with second Piola-Kirchhoff stress tensor $S = \lambda \mathrm{tr}(E)I + 2\mu E$ for
Lagrangian Green strain $E = F^{\top}F - I$ and Lam\'e parameters $\lambda$ and $\mu$, which may be derived from a given
Young's modulus and Poisson ratio.  We solve for displacements $\vu \in \mathcal{V}$, given a constant load vector $b$
imposed on the structure.

\begin{figure}[H]
  \centering
  \includegraphics[width=0.5\textwidth]{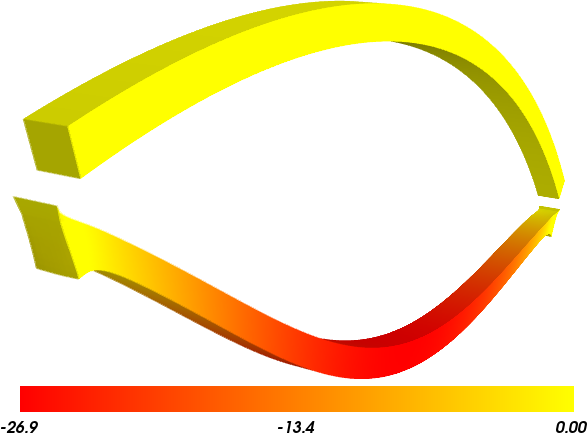}
  \caption{Unstressed and stressed configurations for the elasticity test problem. Coloration indicates vertical displacement in meters.}
  \label{fig:elasticity}
\end{figure}

The domain $\Omega$ is a $60^\circ$ cylindrical arch of inner radius $100$ m.  Homogeneous Dirichlet boundary conditions
are imposed on the outer edge of the ends of the arch.  The goal is to ``snap through'' the arch, causing it to sag
under the load rather than merely compressing it. To cause the sag, we set $b = -\mathbf{e}_y$ and set the Lam\'e
constants consistent with a Young's modulus of 100 and a Poisson ratio of 0.2.  The nonlinearity is highly activated
during the process of snapping through and may be tuned to present a great deal of difficulty to traditional nonlinear
solvers.  The problem is discretized by using hexahedral $\mathbf{Q}_1$ finite elements on a logically structured grid,
deformed to form the arch.  The grid is 401x9x9, and therefore the problem has 97,443 degrees of freedom.  The problem is a
three-dimensional extension of a problem put forth by Wriggers \cite{wriggers2008nonlinear}.  The unstressed and
converged solutions are shown in \figref{fig:elasticity}.
For 3D hexahedral FEM discretizations, \GSN\ and related algorithms are inefficient because for each degree
of freedom each element in the support of the degree of freedom must be visited, resulting in eight visits to an
interior element per sweep.  Instead, we focus on algorithms requiring only a single visit to each cell, restricting us to
function and Jacobian evaluations.  Such approaches are also available in the general unstructured FEM case, and these
experiments may guide users in regimes where no decomposition solvers are available.

For the experiments, we emphasize the role that nonlinear composition and series acceleration may
play in the acceleration of nonlinear solvers. The problem is amenable to \NCG, \ANDERSON, and \NK\ with an
algebraic multigrid preconditioner. We test a number of combinations of these solvers.  Even though we have a
logically structured grid, we approach the problem as if no reasonable grid hierarchy or domain
decomposition were available.  The solver combinations we use are listed in \tabref{tab:elasticitysolvers}. In all
the following tables, Solver denotes outer solver, LPC denotes the linear PC when applicable, NPC denotes the nonlinear
PC when applicable, Side denotes the type of preconditioning, Smooth denotes the level smoothers used in the
multilevel method (\MG/\FAS), and LS denotes line search.
\begin{table}[H]
\centering
\caption{
  Series of solvers for the nonlinear elasticity test problem.}
\label{tab:elasticitysolvers}
\begin{tabular}{r|llllll}
  Name                         & Solver   & LPC     & NPC      & Side & Smooth           & LS \\
  \hline
  $\NCG$                       & \NCG     &         &          &      &                  & \CP \\
  $\NK\pc\MG$                  & \NK      & \MG     & --       & --   & \SOR             & \BT \\
  $\NCG\lp(\NK\pc\MG)$         & \NCG     & --      & $\NK\pc\MG$& L & \SOR             & \CP \\
  $\NGMRES\rp(\NK\pc\MG)$      & \NGMRES  & --      & $\NK\pc\MG$& R & \SOR             & \CP \\
  $\NCG(10)+(\NK\pc\MG)$       & \NCG,\NK & \MG     & --       & --   & \SOR             & \CP/\BT \\
  $\NCG(10)*(\NK\pc\MG)$       & \NCG,\NK & \MG     & --       & --   & \SOR             & \CP/\BT \\
\end{tabular}
\end{table}
For the \NK\ methods, we precondition the inner \GMRES\ solve with a smoothed aggregation algebraic multigrid method
provided by the GAMG package in PETSc.  The relative tolerance for the \GMRES\ solve is $10^{-3}$.  For all
instances of the \NCG\ solver, the \CP\ line search initial guess for $\lambda$ is the final value from the previous
nonlinear iteration.  In the case of $\NGMRES\rp(\NK\pc\MG)$ the inner line search is \LT.  \NGMRES\ stores up to 30
previous solutions and residuals for all experiments.  The outer line search for $\NCG\lp(\NK\pc\MG)$ is a second-order
secant approximation rather than first-order as depicted in \agmref{alg:cplinesearch}.  For the composite examples, 10
iterations of $\NCG$ are used as one of the subsolvers, denoted $\NCG(10)$.  The additive \composition's
weights are determined by the \ANDERSON\ minimization mechanism as described in \secref{sec:anderson}.

The example, \href{http://www.mcs.anl.gov/petsc/petsc-current/src/snes/examples/tutorials/ex16.c.html}{SNES ex16.c},
can be run directly by using a default PETSc installation.  The command line used for these experiments is

\begin{framed}
\begin{alltt}\scriptsize
./ex16 -da_grid_x 401 -da_grid_y 9 -da_grid_z 9 -height 3 -width 3
  -rad 100 -young 100 -poisson 0.2 -loading -1 -ploading 0
\end{alltt}
\end{framed}

\begin{figure}[H]
  \centering
  \includegraphics[width=0.5\textwidth]{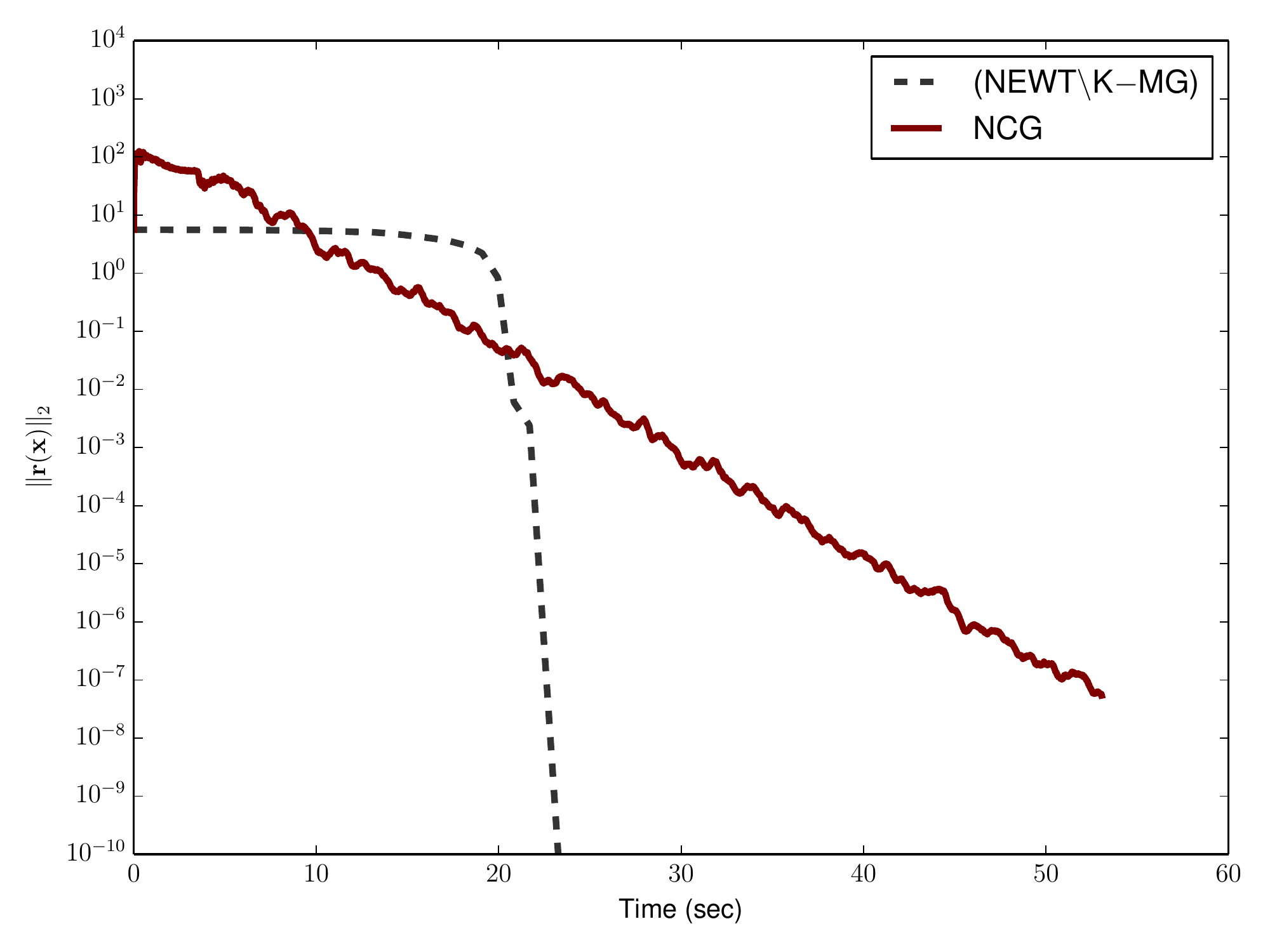}
  \caption{$(\NK\pc\MG)$ and \NCG\ convergence.}
  \label{fig:ex16Unpre}
\end{figure}
\begin{table}[H]
\centering
  \caption{$(\NK\pc\MG)$ and \NCG\ results.}
\label{tab:elasticity:np}
\begin{tabular}{r|lllllll}
Solver & T & N. It & L. It & Func & Jac & PC & NPC \\
\hline
$(\NK\pc\MG)$ & 23.43 & 27 & 1556 & 91 & 27 & 1618 & -- \\
$\NCG$ & 53.05 & 4495 & 0 & 8991 & -- & -- & -- \\
\end{tabular}
\end{table}
The solvers shown in \tabref{tab:elasticity:np} converge slowly.  $\NK\pc\MG$ takes 27 nonlinear iterations and 1,618
multigrid V-cycles to reach convergence.  \NCG\ requires 8,991 function evaluations and since it is unpreconditioned,
scales unacceptably with problem size with respect to number of iterations.  Note in \figref{fig:ex16Unpre} that while
\NCG\ takes an initial jump and then decreases at a constant rate, $\NK\pc\MG$ is near stagnation until the end,
when it suddenly drops into the basin of attraction.
\begin{figure}[H]
    \centering
    \includegraphics[width=0.5\textwidth]{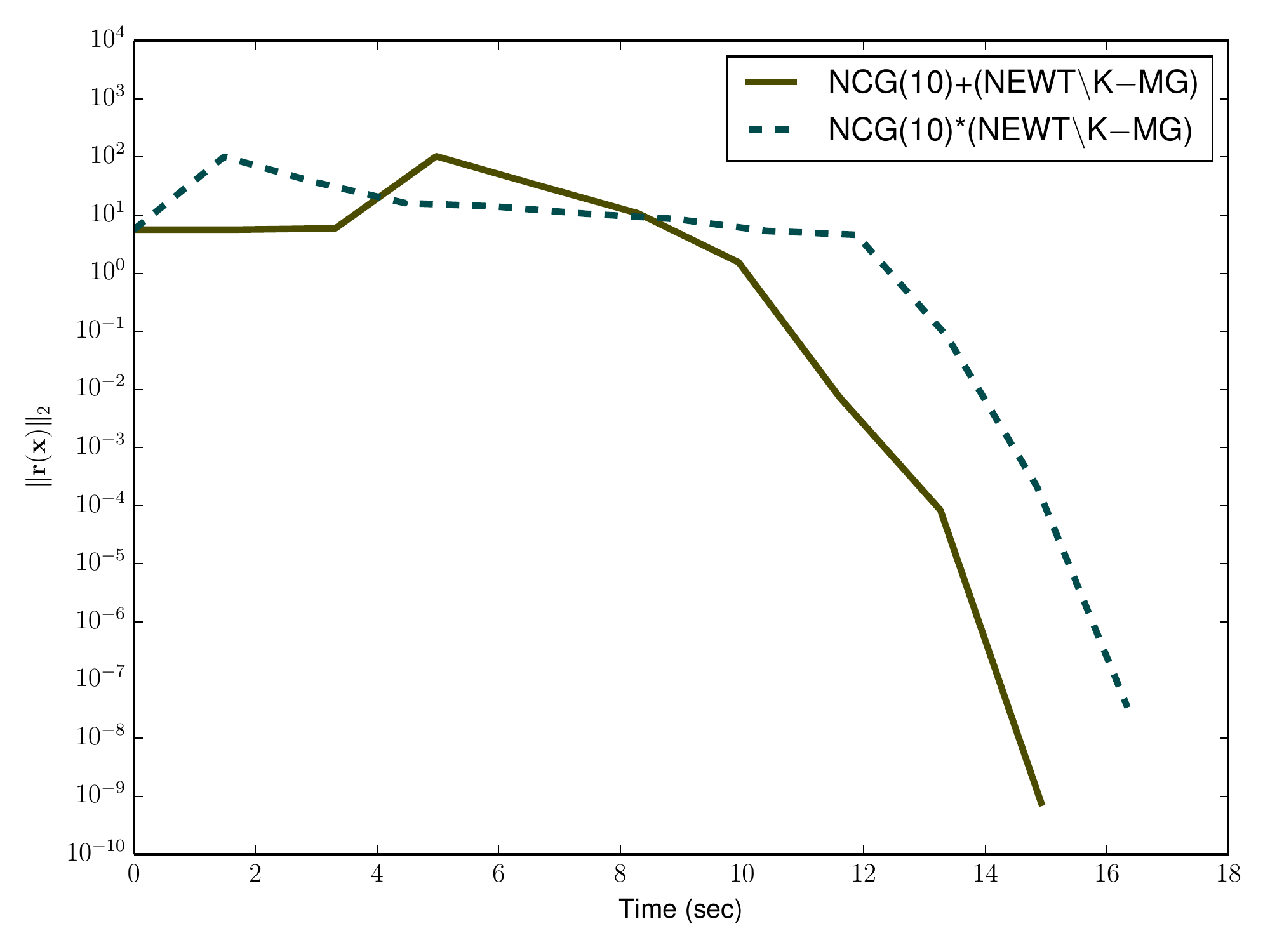}
    \caption{$\NCG(10)+(\NK\pc\MG)$ and $\NCG(10)*(\NK\pc\MG)$  convergence.}
    \label{fig:ex16Composed}
\end{figure}
\begin{table}[H]
\centering
\caption{$\NCG(10)+(\NK\pc\MG)$ and $\NCG(10)*(\NK\pc\MG)$  results.}
\label{tab:elasticity:composed}
\begin{tabular}{r|lllllll}
Solver & T & N. It & L. It & Func & Jac & PC & NPC \\
\hline
$\NCG(10)+(\NK\pc\MG)$ & 14.92 & 9 & 459 & 218 & 9 & 479 & -- \\
$\NCG(10)*(\NK\pc\MG)$ & 16.34 & 11 & 458 & 251 & 11 & 477 & -- \\
\end{tabular}
\end{table}
Results for \composition are listed in \tabref{tab:elasticity:composed}.  Additive and multiplicative composite
combination provides roughly similar speedups for the problem, with more total outer iterations in the additive case.
As shown in \figref{fig:ex16Composed}, neither the additive nor multiplicative methods stagnate.  After the
same initial jump and convergence at the pace of \NCG, the basin of attraction is reached, and the entire iteration
converges quadratically as should be expected with \NK.

\begin{figure}[H]
    \centering
    \includegraphics[width=0.5\textwidth]{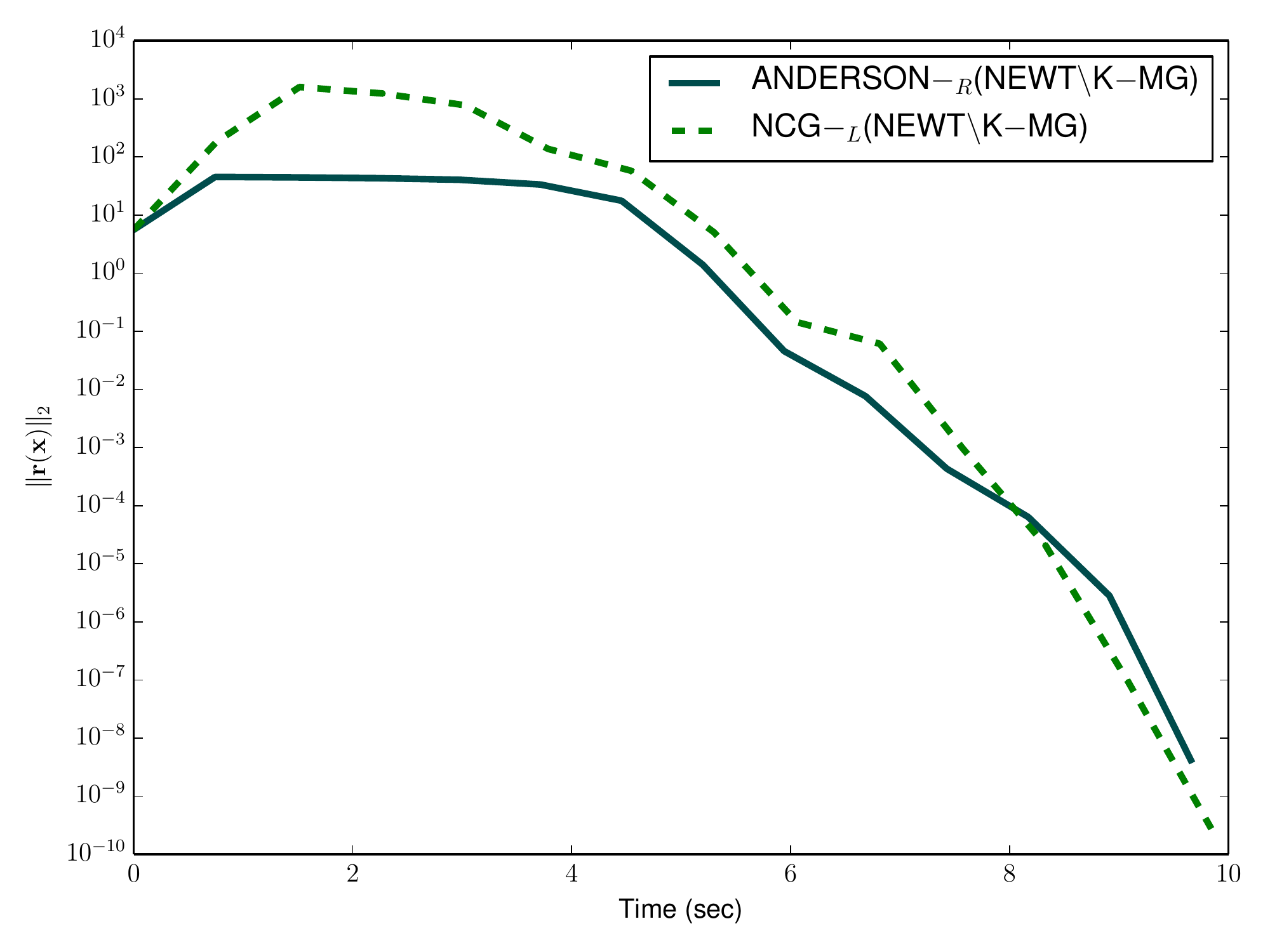}
    \caption{$\NGMRES\rp(\NK\pc\MG)$ and $\NCG\lp(\NK\pc\MG)$ convergence.}
    \label{fig:ex16Pre}
\end{figure}
\begin{table}[H]
\centering
\caption{$\NGMRES\rp(\NK\pc\MG)$ and $\NCG\lp(\NK\pc\MG)$ results.}
\label{tab:elasticity:pc}
\begin{tabular}{r|lllllll}
Solver & T & N. It & L. It & Func & Jac & PC & NPC \\
\hline
$\NGMRES\rp(\NK\pc\MG)$ & 9.65 & 13 & 523 & 53 & 13 & 548 & 13 \\
$\NCG\lp(\NK\pc\MG)$ & 9.84 & 13 & 529 & 53 & 13 & 554 & 13 \\
\end{tabular}
\end{table}
Left-preconditioning of \NCG\ and \ANDERSON\ with \NK\pc\MG\ provides even greater benefits, as shown in
\figref{fig:ex16Pre} and \tabref{tab:elasticity:pc}.  The number of iterations is nearly halved from the
unpreconditioned case, and the number of linear preconditioner applications is only slightly increased compared with the
\composition cases.  In \figref{fig:ex16Pre}, we see that the Newton's method-preconditioned nonlinear Krylov solvers
never dive as \NK\ did but maintain rapid, mostly constant convergence instead.  This constant convergence is a
significantly more efficient path to solution than the $\NK\pc\MG$ solver alone.
\subsection{Driven Cavity}
\label{sec:liddriven}
The driven cavity formulation used for the next set of experiments can be stated as
\begin{eqnarray}
\label{eq:velocityvorticity}
  - \Delta \vu    -             \nabla \times \Omega                               &=& 0  \\
  - \Delta \Omega +             \nabla \cdot (\vu \Omega) - \mathrm{Gr} \nabla_x T &=& 0  \\
  - \Delta T      + \mathrm{Pr} \nabla \cdot \vu                                   &=& 0.
\end{eqnarray}
The fluid motion is driven in part by a moving lid and in part by buoyancy.  A Grashof number
$Gr=2e4$, Prandtl number $Pr=1$, and lid velocity of 100 are used in these experiments.  No-slip, rigid-wall Dirichlet
conditions are imposed for $\vu$. Dirichlet conditions are used for $\Omega$, based on the definition of vorticity,
$\Omega = \nabla \times \vu$, where along each constant coordinate boundary the tangential derivative is
zero. Dirichlet conditions are used for $T$ on the left and right walls, and insulating homogeneous Neumann conditions
are used on the top and bottom walls. A finite-difference approximation with the usual 5-point stencil is used to
discretize the boundary value problem in order to obtain a nonlinear system of equations.  Upwinding is used for the divergence
(convective) terms and central differencing for the gradient (source) terms.

In these experiments, we emphasize the use of \ANDERSON, \MG\ or \FAS, and \GSN.  The solver combinations are listed in
\tabref{tab:fassolvers}.  The multilevel methods use a simple series of structured grids, with the smallest being 5x5
and the largest being 257x257.  With four fields per grid point, we have an overall system size of 264,196
unknowns.  In $\NK-\MG$, the Jacobian is rediscretized on each level.  \GSN\ is used as the smoother for \FAS\ and
solves the four-component block problem corresponding to \eqnref{eq:velocityvorticity} per grid point in a simple sweep
through the processor local part of the domain.

The example, \href{http://www.mcs.anl.gov/petsc/petsc-current/src/snes/examples/tutorials/ex19.c.html}{SNES ex19.c},
can be run directly by using a default PETSc installation.  The command line used for these experiments is

\begin{figure}[H]
\centering
\begin{subfigure}{0.3\textwidth}
    \label{fig:lidxvelocity}
    \includegraphics[width=\textwidth]{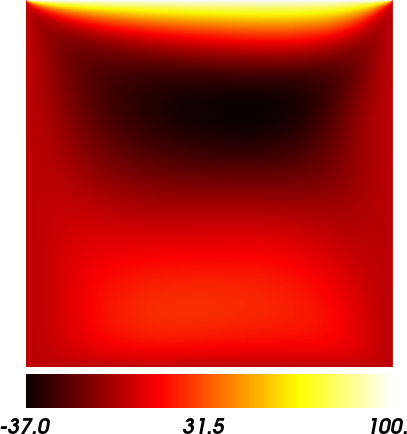}
\end{subfigure}
\begin{subfigure}{0.3\textwidth}
    \label{fig:lidyvelocity}
    \includegraphics[width=\textwidth]{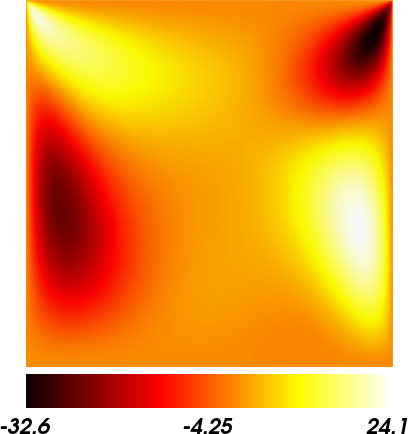}
\end{subfigure}
\begin{subfigure}{0.3\textwidth}
    \label{fig:lidtemperature}
    \includegraphics[width=\textwidth]{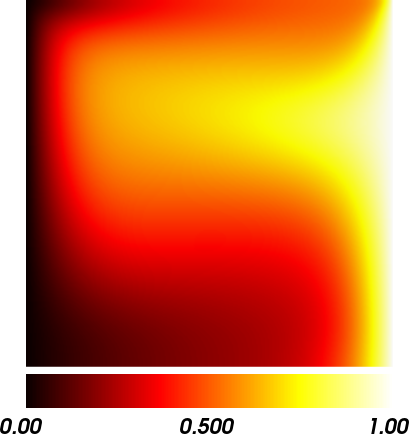}
\end{subfigure}
\caption{$\vu_x$, $\vu_y$, and temperature profiles for the lid driven cavity solution.}
\label{fig:lid}
\end{figure}

\begin{framed}
\begin{alltt}\scriptsize
  ./ex19 -da_refine 6 -da_grid_x 5 -da_grid_y 5 -grashof 2e4 -lidvelocity 100 -prandtl 1.0
\end{alltt}
\end{framed}
\noindent and the converged solution using these options is shown in \figref{fig:lid}.
\begin{table}[H]
\centering
\caption{
  Solvers for the lid driven cavity problem.
}
\label{tab:fassolvers}
\begin{tabular}{r|llllll}
  Name                        & Solver  & LPC & NPC  & Side & \MG/\FAS Smooth. & LS  \\
  \hline
  $(\NK\pc\MG)$               & \NK     & \MG     & --       & --   & \SOR             & \BT \\
  $\NGMRES\rp(\NK\pc\MG$)     & \NGMRES & \MG     & \NK      & R    & \SOR             & --  \\
  \FAS                        & \FAS    & --      & --       & --   & \GSN             & --  \\
  $\NRICH\lp\FAS$             & \NRICH  & --      & \FAS     & L    & \GSN             & \LT \\
  $\NGMRES\rp\FAS$            & \NGMRES & --      & \FAS     & R    & \GSN             & --  \\
  $\FAS*(\NK\pc\MG)$          & \FAS/\NK& \MG     & --       & --   & \SOR/\GSN        & \BT \\
  $\FAS+(\NK\pc\MG)$          & \FAS/\NK& \MG     & --       & --   & \SOR/\GSN        & \BT \\
\end{tabular}
\end{table}
All linearized problems in the \NK\ iteration are solved by using GMRES with geometric multigrid preconditioning to a
relative tolerance of $10^{-8}$.  There are five levels, with Chebychev\pc\SOR\ smoothers on each level.  For \FAS, the
smoother is \GSN(5).  The coarse-level smoother is five iterations of \NK-LU; chosen for robustness.  In the case of
$\NGMRES\rp(\NK\pc\MG)$ and $\FAS+(\NK\pc\MG)$, the $\NK\pc\MG$ step is damped by one half.  Note that acting alone,
this damping would cut the rate of convergence to linear with a rate constant $\frac{1}{2}$.  The ideal step size is
recovered by the minimization procedure.  $\NK\pc\MG$ is undamped in all other tests.
In these experiments we emphasize the total number of V-cycles, linear and nonlinear, since they
dominate the runtime and contain the lion's share of the communication and floating-point operations.
\begin{figure}[H]
    \centering
    \includegraphics[width=0.5\textwidth]{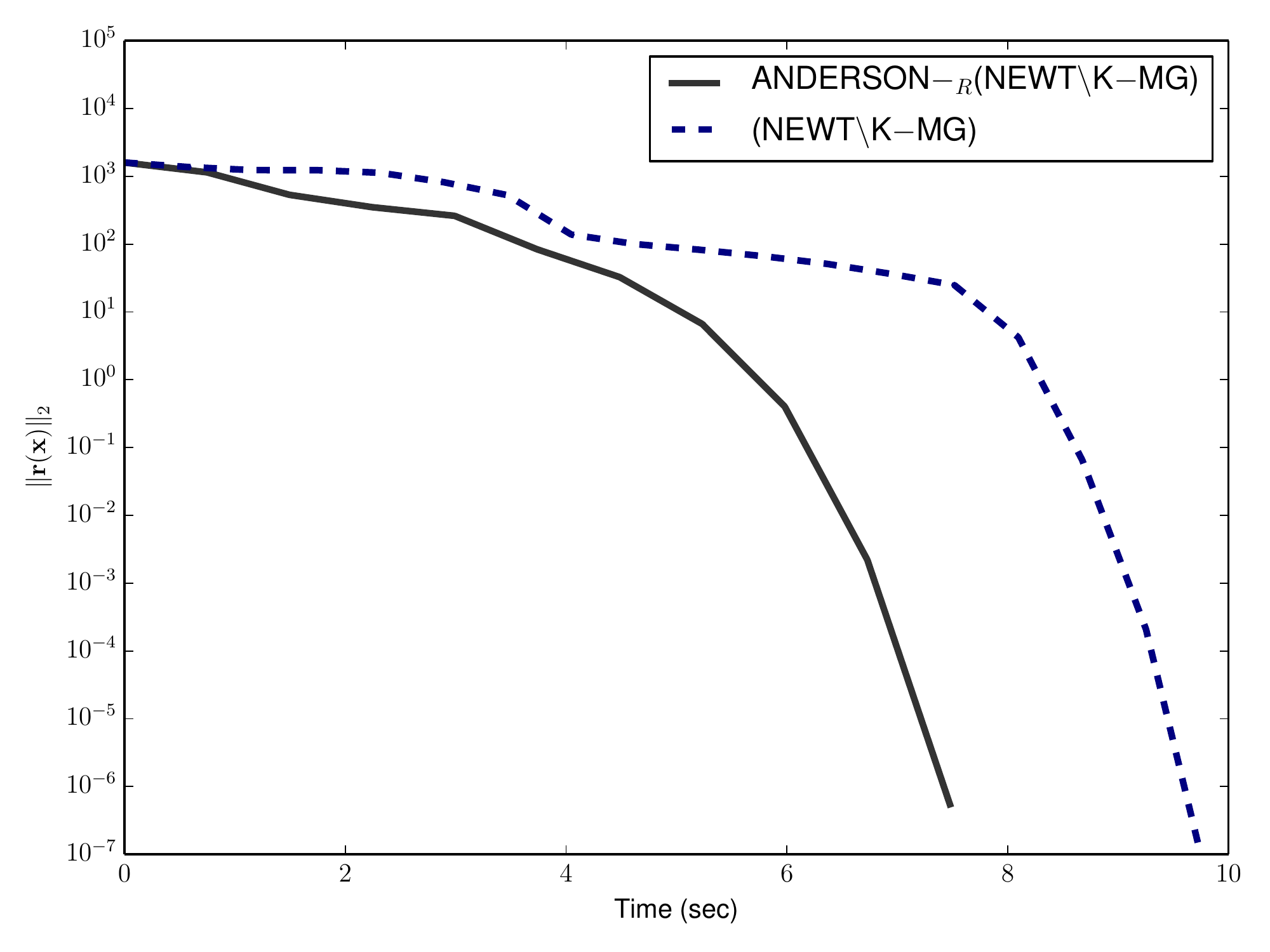}
\caption{$\NGMRES\rp(\NK\pc\MG)$ and $(\NK\pc\MG)$ convergence.}
    \label{fig:ex19Newton}
\end{figure}
\begin{table}[H]
\centering
\caption{$\NGMRES\rp(\NK\pc\MG)$ and $(\NK\pc\MG)$ results.}
\label{tab:bouyant:unpre}
\begin{tabular}{r|lllllll}
Solver & T & N. It & L. It & Func & Jac & PC & NPC \\
\hline
$\NGMRES\rp(\NK\pc\MG)$ & 7.48 & 10 & 220 & 21 & 50 & 231 & 10 \\
$(\NK\pc\MG)$ & 9.83 & 17 & 352 & 34 & 85 & 370 & -- \\
\end{tabular}
\end{table}
In \tabref{tab:bouyant:unpre}, we see that Newton's method converges in 17 iterations, with 370 V-cycles.  Using
\NGMRES\ instead of a line search provides some benefit, as $\NGMRES\pc(\NK\pc\MG)$ takes 231 V-cycles and 10 iterations.
\begin{figure}[H]
    \centering
    \includegraphics[width=0.5\textwidth]{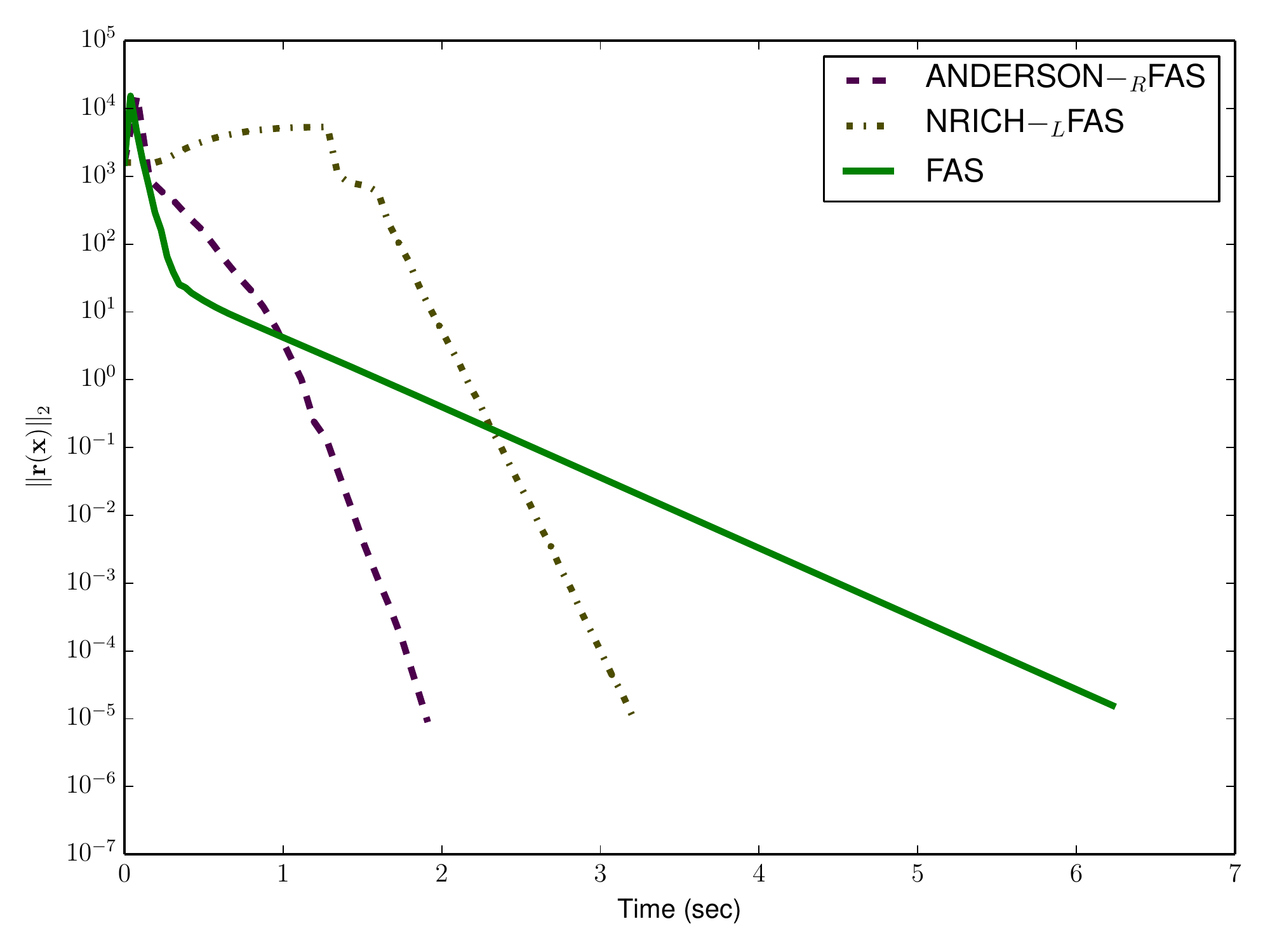}
\caption{$\NGMRES\rp\FAS$, $\NRICH\lp\FAS$, and $\FAS$ convergence.}
    \label{fig:ex19FAS}
\end{figure}
\begin{table}[H]
\centering
\caption{$\NGMRES\rp\FAS$, $\NRICH\lp\FAS$, and $\FAS$ results.}
\label{tab:bouyant:fas}
\begin{tabular}{r|lllllll}
Solver & T & N. It & L. It & Func & Jac & PC & NPC \\
\hline
$\NGMRES\rp\FAS$ & 1.91 & 24 & 0 & 447 & 83 & 166 & 24 \\
$\NRICH\lp\FAS$ & 3.20 & 50 & 0 & 1180 & 192 & 384 & 50 \\
$\FAS$ & 6.23 & 162 & 0 & 2382 & 377 & 754 & -- \\
\end{tabular}
\end{table}

$\NRICH\lp\FAS$ takes 50 V-cycles, at the expense of
three more fine-level function evaluations per iteration.  $\NGMRES\rp\FAS$ reduces the number of V-cycles to 24 at the expense
of more communication.
\begin{figure}[H]
    \centering
    \includegraphics[width=0.5\textwidth]{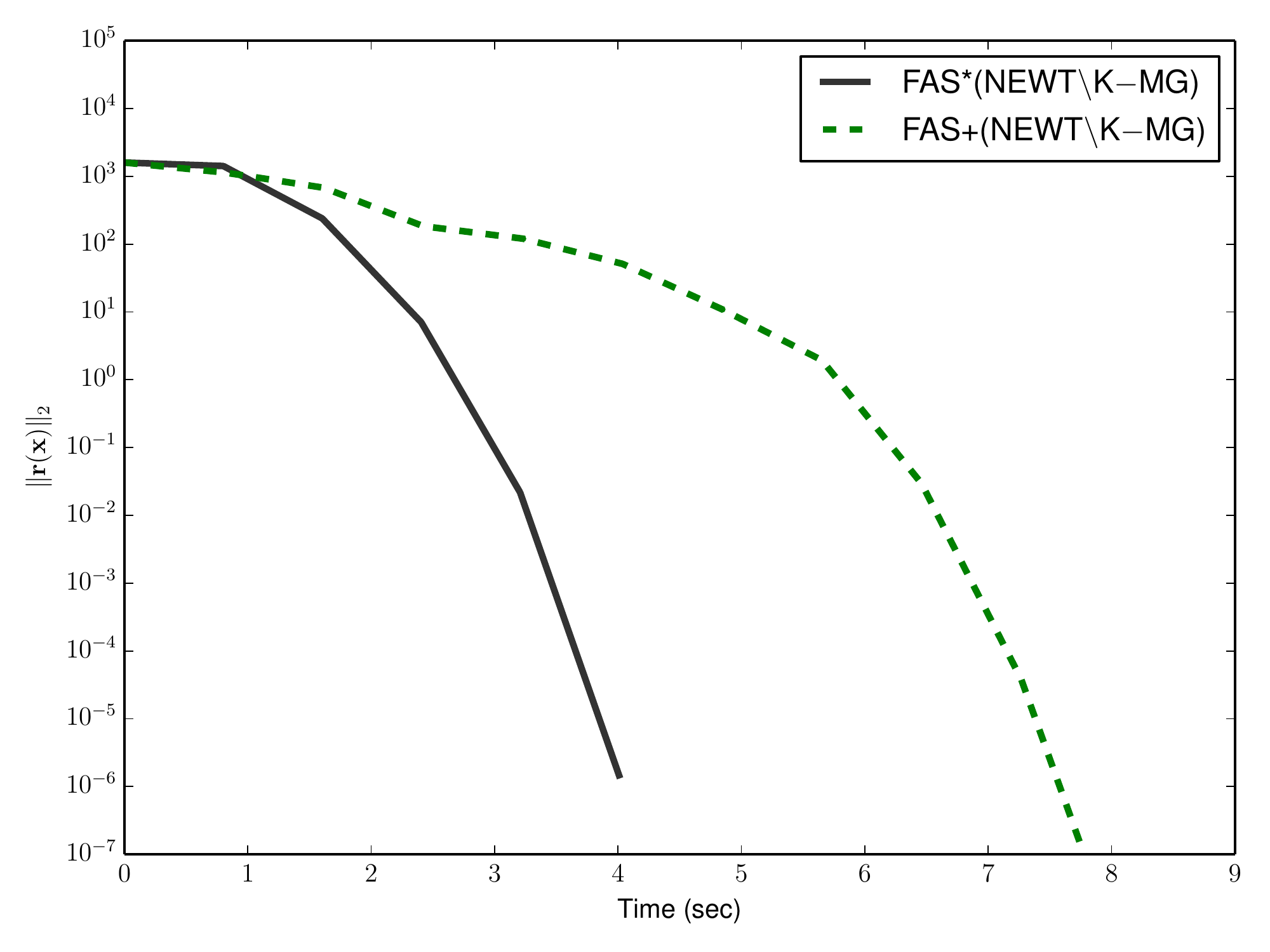}
    \caption{$\FAS*(\NK\pc\MG)$ and $\FAS+(\NK\pc\MG)$ convergence.}
    \label{fig:ex19Composed}
\end{figure}
\begin{table}[H]
\centering
\caption{$\FAS*(\NK\pc\MG)$ and $\FAS+(\NK\pc\MG)$ results.}
\label{tab:bouyant:composed}
\begin{tabular}{r|lllllll}
Solver & T & N. It & L. It & Func & Jac & PC & NPC \\
\hline
$\FAS*(\NK\pc\MG)$ & 4.01 & 5 & 80 & 103 & 45 & 125 & -- \\
$\FAS+(\NK\pc\MG)$ & 8.07 & 10 & 197 & 232 & 90 & 288 & -- \\
\end{tabular}
\end{table}
Composed nonlinear methods are shown in \tabref{tab:bouyant:composed}.  Multiplicative \composition consisting
of \FAS\ and \NK\pc\MG\ reduces the total number of V-cycles to 130.  Additive \composition using the
least-squares minimization is less effective, taking 298 V-cycles.  Note in \figref{fig:ex19Composed} that both the
additive and multiplicative solvers show that combining \FAS\ and \NK-\MG\ may speed solution, with
multiplicative combination being significantly more effective.

\subsection{Tuning the Solvers to Obtain Convergence}
We now show how the composed and preconditioned solves may be tuned for more difficult nonlinearities where the basic methods fail to converge using the same model problem.
 For Grashof number $\mathrm{Gr} < 10^4$ and Prandtl number $\mathrm{Pr} =
1.0$, Newton's method converges well:
\begin{framed}
\begin{alltt}\scriptsize
lid velocity = 100, prandtl # = 1, grashof # = 10000
  0 SNES Function norm 715.271 
  1 SNES Function norm 623.41 
  2 SNES Function norm 510.225 
   \vdots
  6 SNES Function norm 0.269179 
  7 SNES Function norm 0.00110921 
  8 SNES Function norm 1.12763e-09 
Nonlinear solve converged due to CONVERGED_FNORM_RELATIVE iterations 8
Number of SNES iterations = 8
\end{alltt}
\end{framed}
For higher Grashof number, Newton's method stagnates
\begin{framed}
\begin{alltt}\scriptsize
./ex19 -lidvelocity 100 -grashof 5e4 -da_refine 4 -pc_type lu -snes_monitor_short -snes_converged_reason
lid velocity = 100, prandtl # = 1, grashof # = 50000
  0 SNES Function norm 1228.95 
  1 SNES Function norm 1132.29 
            \vdots
 29 SNES Function norm 580.937 
 30 SNES Function norm 580.899 
\end{alltt}
\end{framed}
It also fails with the standard continuation strategy from coarser meshes (not shown).
We next try nonlinear multigrid, in the hope that multiple updates from a coarse solution are
sufficient, using \GS\ as the nonlinear smoother,
\begin{framed}
\begin{alltt}\scriptsize
./ex19 -lidvelocity 100 -grashof 5e4 -da_refine 4 -snes_monitor_short -snes_converged_reason -snes_type fas 
       -fas_levels_snes_type ngs -fas_levels_snes_max_it 6 -snes_max_it 25
lid velocity = 100, prandtl # = 1, grashof # = 50000
  0 SNES Function norm 1228.95 
  1 SNES Function norm 574.793 
  2 SNES Function norm 513.02 
  3 SNES Function norm 216.721 
  4 SNES Function norm 85.949 
  5 SNES Function norm 108.24 
  6 SNES Function norm 207.469 
  \vdots
 22 SNES Function norm 131.866 
 23 SNES Function norm 114.817 
 24 SNES Function norm 71.4699 
 25 SNES Function norm 63.5413 
\end{alltt}
\end{framed}
\noindent But the residual norm just jumps around and the method does not converge. We then accelerate the method with \ANDERSON\ and obtain convergence.
\begin{framed}
\begin{alltt}\scriptsize
./ex19 -lidvelocity 100 -grashof 5e4 -da_refine 4 -snes_monitor_short -snes_converged_reason -snes_type anderson
       -npc_snes_max_it 1 -npc_snes_type fas  -npc_fas_levels_snes_type ngs -npc_fas_levels_snes_max_it 6
lid velocity = 100, prandtl # = 1, grashof # = 50000
  0 SNES Function norm 1228.95 
  1 SNES Function norm 574.793 
  2 SNES Function norm 345.592 
  3 SNES Function norm 155.476 
  4 SNES Function norm 70.2302 
  5 SNES Function norm 40.3618 
  6 SNES Function norm 29.3065 
  7 SNES Function norm 14.2497 
  8 SNES Function norm 4.80462 
  9 SNES Function norm 4.15985 
 10 SNES Function norm 2.13428 
 11 SNES Function norm 1.57717 
 12 SNES Function norm 0.60919 
 13 SNES Function norm 0.150496 
 14 SNES Function norm 0.0355709 
 15 SNES Function norm 0.00705481 
 16 SNES Function norm 0.00164509 
 17 SNES Function norm 0.000464835 
 18 SNES Function norm 6.02035e-05 
 19 SNES Function norm 1.11713e-05 
Nonlinear solve converged due to CONVERGED_FNORM_RELATIVE iterations 19
Number of SNES iterations = 19
\end{alltt}
\end{framed}
We can restore the convergence to a few iterates by increasing the power of the nonlinear smoothers in \FAS\. We replace \GSN\ by six iterations of Newton using the default linear solver of GMRES plus ILU(0). Note that as a solver alone this does not converge but it performs very well as a smoother for \FAS.
\begin{framed}
\begin{alltt}\scriptsize
./ex19 -lidvelocity 100 -grashof 5e4 -da_refine 4 -snes_monitor_short -snes_converged_reason -snes_type anderson
       -npc_snes_max_it 1 -npc_snes_type fas -npc_fas_levels_snes_type newtonls -npc_fas_levels_snes_max_it 6 
       -npc_fas_levels_snes_linesearch_type basic -npc_fas_levels_snes_max_linear_solve_fail 30 
       -npc_fas_levels_ksp_max_it 20
lid velocity = 100, prandtl # = 1, grashof # = 50000
  0 SNES Function norm 1228.95 
  1 SNES Function norm 0.187669 
  2 SNES Function norm 0.0319743 
  3 SNES Function norm 0.00386815 
  4 SNES Function norm 2.24093e-05 
  5 SNES Function norm 5.38246e-08 
Nonlinear solve converged due to CONVERGED_FNORM_RELATIVE iterations 5
Number of SNES iterations = 5
\end{alltt}
\end{framed}
Thus we have demonstrated how one may experimentally add composed solvers to go from complete lack of convergence to convergence with a small number of iterations.

\subsection{$p$-Laplacian}
\label{sec:plaplacian}
\begin{figure}[H]
  \centering
  \begin{subfigure}{0.3\textwidth}
    \label{fig:plaplacianinitial}
    \includegraphics[width=\textwidth]{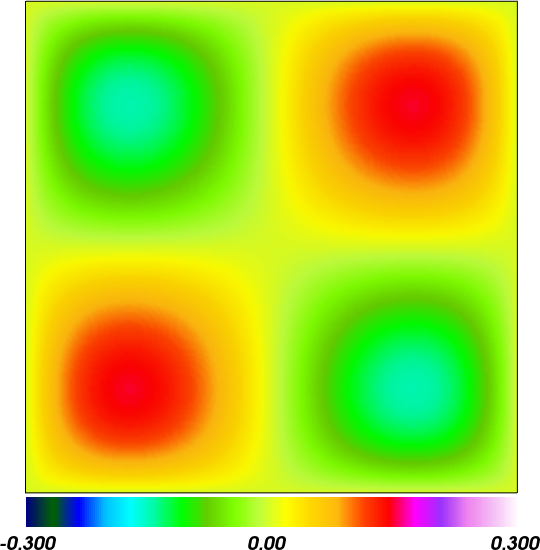}
  \end{subfigure}
  \begin{subfigure}{0.3\textwidth}
    \label{fig:plaplacianfinal}
    \includegraphics[width=\textwidth]{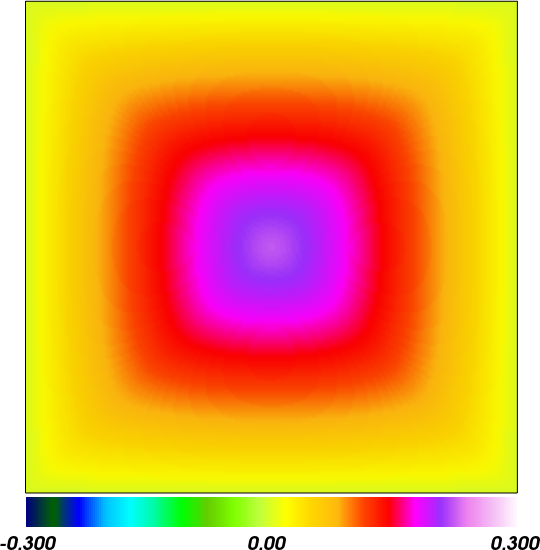}
  \end{subfigure}
  \caption{Initial and converged solutions to the $p=5$ $p$-Laplacian.}
  \label{fig:plaplacian}
\end{figure}

The regularized $p$-Laplacian formulation used for these experiments is
\begin{equation*}
        -\nabla\cdot ((\epsilon^2 + \frac{1}{2}|\nabla u|^2)^{(p-2)/2} \nabla u) = c,
\end{equation*}
\noindent where $\epsilon = 10^{-5}$ is the regularization parameter and $p$ the exponent of the Laplacian. When $p = 2$
the $p$-Laplacian reduces to the Poisson equation. We consider the case where $p = 5$ and $c = 0.1$.  The domain is
$[-1, 1]\times[-1, 1]$ and the initial guess is $u_0(x,y) = xy(1 - x^2)(1 - y^2)$.  The grid used is 385x385, leading to
a total of 148,225 unknowns in the system.  The initial and converged solutions are shown in \figref{fig:plaplacian}.
The example, \href{http://www.mcs.anl.gov/petsc/petsc-current/src/snes/examples/tutorials/ex15.c.html}{SNES ex15.c},
can be run directly by using a default PETSc installation.  The command line used for these experiments is
\begin{framed}
\begin{alltt}\scriptsize
./ex15 -da_refine 7 -da_overlap 6 -p 5.0 -lambda 0.0 -source 0.1
\end{alltt}
\end{framed}
\begin{table}[H]
\centering
\caption{Solvers for the $p$-Laplacian problem.
  SubPC denotes the linear solver at the block level.}
\label{tab:plaplacian:solvers}
\begin{tabular}{r|llllll}
  Name                      & Solver & LPC  & NPC  & Side & SubPC   & LS \\
  \hline
  $\NK\pc\ASM$              & \NK      & \ASM    & --       & --   & \LU             & \BT \\
  $\QN$                     & \QN      & --      & --       & --   &                 & \CP \\
  $\RAS$                    & \RAS     & --      & --       & --   & \LU             & \CP \\
$\RAS+(\NK\pc\ASM)$           & \RAS/\NK & \ASM    & --       & --   & \LU             & \BT \\
$\RAS*(\NK\pc\ASM)$           & \RAS/\NK & \ASM    & --       & --   & \LU             & \BT \\
  $\ASPIN$                  & \NK      & --      & \NASM    & L    & \LU             & \BT \\
$\NRICH\lp(\RAS)$           & \NRICH   & --      & \RAS     & L    & \LU             & \CP \\
$\QN\lp(\RAS)$              & \QN      & --      & \RAS     & L    & \LU             & \CP \\
\end{tabular}
\end{table}
We concentrate on additive Schwarz and \QN\ methods, as listed in \tabref{tab:plaplacian:solvers}. We will show that 
combination has distinct advantages.  This problem has difficult local nonlinearity that may impede global
solvers, so the local solvers should be an efficient remedy.  We also consider a \composition of Newton's
method with \RAS\ and nonlinear preconditioning of Newton's method in the form of \ASPIN.  \NASM, \RAS, and linear
\ASM\ all have single subdomains on each processor that overlap each other by six grid points.  We use the L-BFGS variant
of \QN\ as explained in \secref{sec:qn}.  Note that the function evaluation is inexpensive for this problem and
methods that use many evaluations perform well.  The advantage of such methods is exacerbated by the large number of
iterations required by Newton's method.  All the \NK\ solvers use an inner \GMRES\ iteration with a relative tolerance
of $10^{-5}$.  \GMRES\ for \ASPIN\ is set to have an inner tolerance of $10^{-3}$.
\begin{figure}[H]
    \centering
    \includegraphics[width=0.5\textwidth]{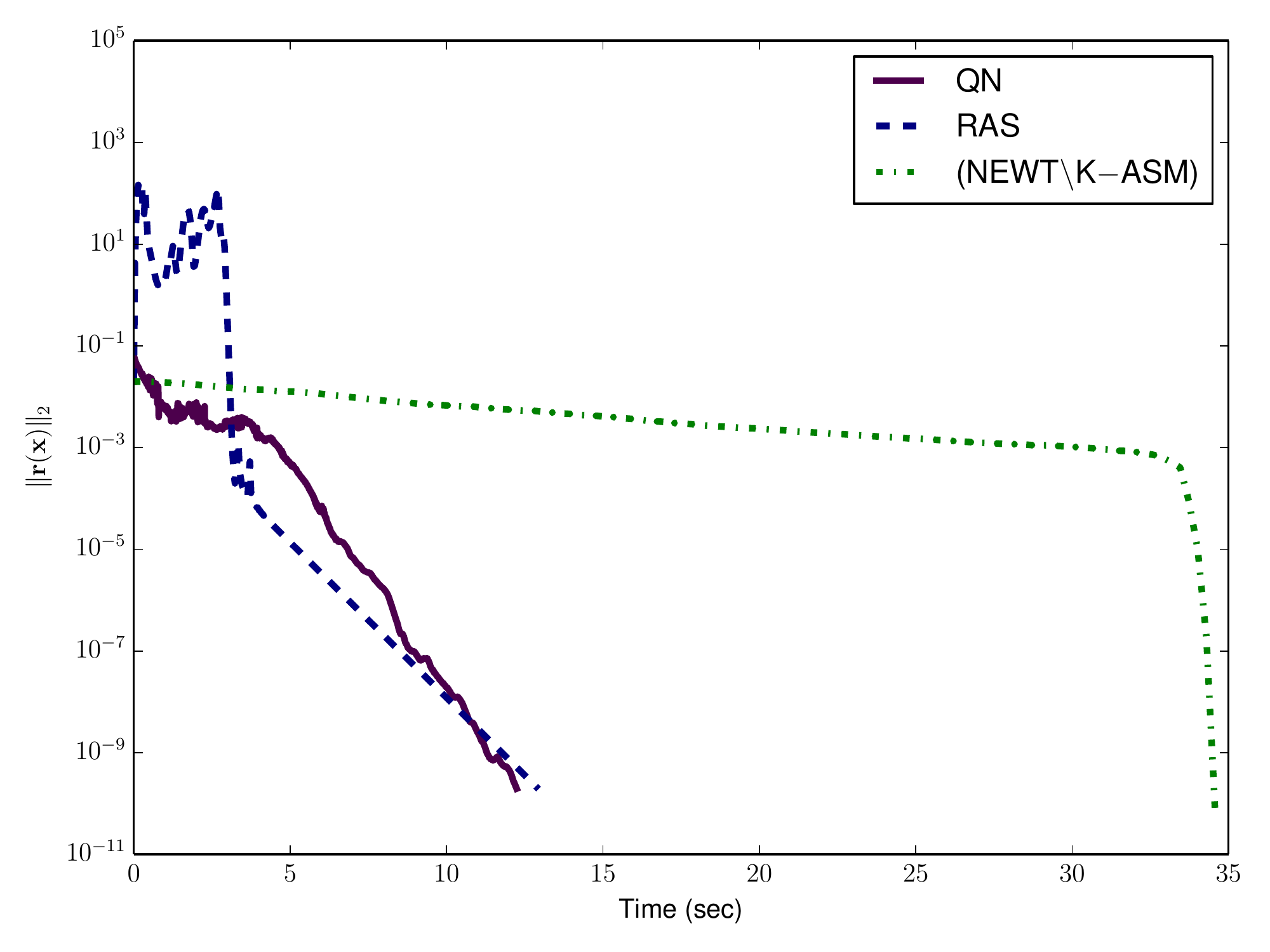}
    \caption{\QN, \RAS, and $(\NK\pc\ASM)$ convergence.}
    \label{fig:ex15Unpre}
\end{figure}
\begin{table}[H]
  \centering
  \caption{\QN, \RAS, and $(\NK\pc\ASM)$ results.}
  \label{tab:plaplacian:np}
\begin{tabular}{r|lllllll}
Solver & T & N. It & L. It & Func & Jac & PC & NPC \\
\hline
$\QN$          & 12.24 & 2960 & 0    & 5921 & --  & --   & -- \\
$\RAS$         & 12.94 & 352  & 0    & 1090 & 352 & 352  & -- \\
$(\NK\pc\ASM)$ & 34.57 & 124  & 3447 & 423  & 124 & 3574 & -- \\
\end{tabular}
\end{table}
\tabref{tab:plaplacian:np} contains results for the uncomposed solvers.  The default solver, $\NK\pc\ASM$, takes a large
number of outer Newton iterations and inner $\GMRES-\ASM$ iterations.  In addition, the line search is
consistently activated, causing an average of more than three function evaluations per iteration.  Unpreconditioned
\QN\ is able to converge to the solution efficiently and will be hard to beat.  \QN\ stores up to 10 previous solutions
and residuals.  \RAS\ set to do one local Newton iteration per subdomain per outer iteration proves to be
efficient on its own after some initial problems, as shown in \figref{fig:ex15Unpre}.
\begin{figure}[H]
    \centering
    \includegraphics[width=0.5\textwidth]{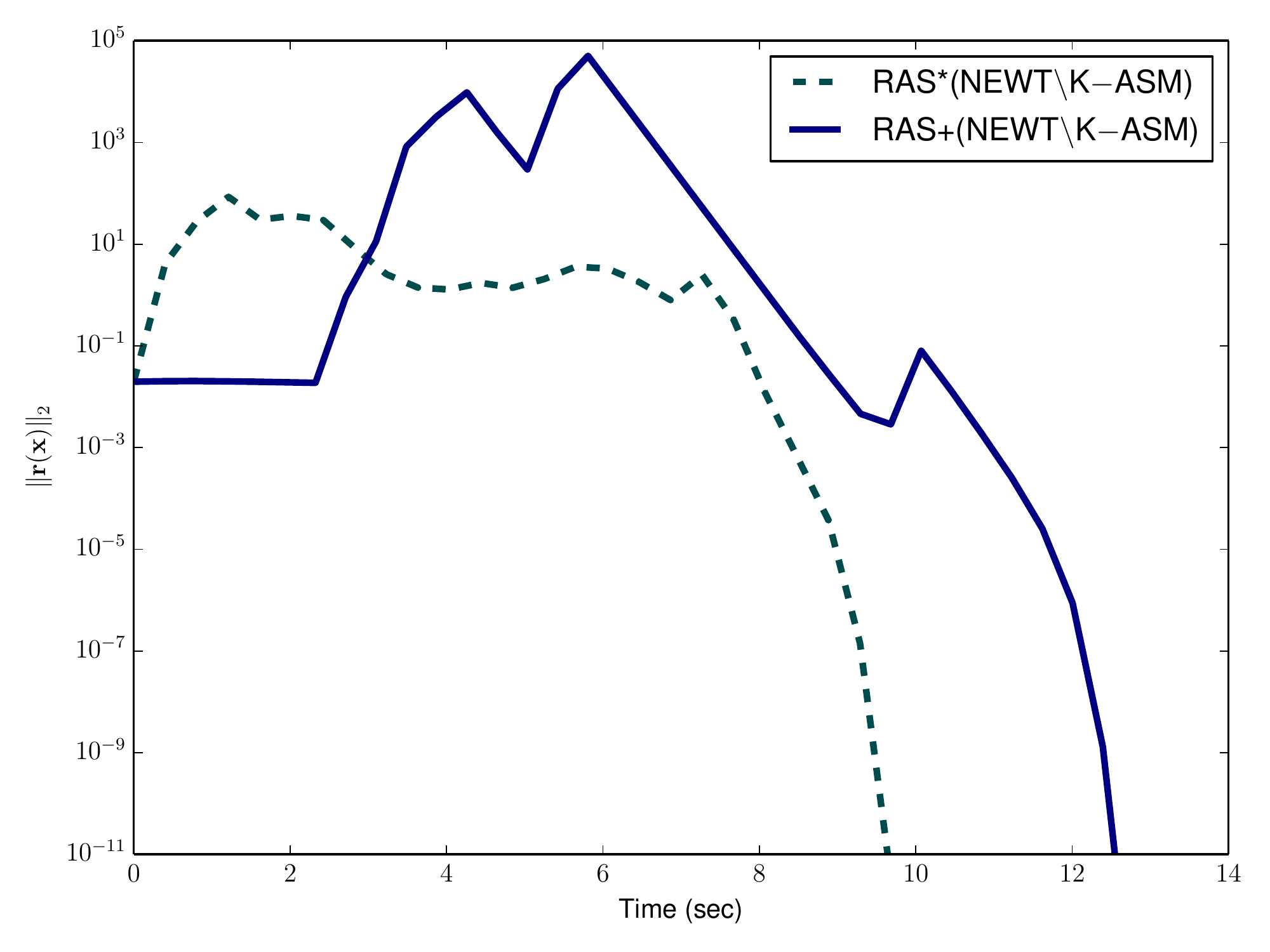}
    \caption{$\RAS*(\NK\pc\ASM)$ and $\RAS+(\NK\pc\ASM)$ convergence.}
    \label{fig:ex15Newton}
\end{figure}
\begin{table}[H]
  \centering
  \caption{$\RAS*(\NK\pc\ASM)$ and $\RAS+(\NK\pc\ASM)$ results.}
\label{tab:plaplacian:composed}
\begin{tabular}{r|lllllll}
Solver & T & N. It & L. It & Func & Jac & PC & NPC \\
\hline
$\RAS*(\NK\pc\ASM)$ & 9.69  & 24 & 750 & 142 & 48 & 811  & -- \\
$\RAS+(\NK\pc\ASM)$ & 12.78 & 33 & 951 & 232 & 66 & 1023 & -- \\
\end{tabular}
\end{table}
In \tabref{tab:plaplacian:composed} and \figref{fig:ex15Newton}, we show how Newton's method may be dramatically improved
by \composition with \RAS\ and \NASM.  The additive \composition reduces the number of outer
iterations substantially.  The multiplicative \composition is even more effective, reducing the number of outer
iterations by a factor of 5 and decreasing runtime by a factor of around 4.
\begin{figure}[H]
    \centering
  \includegraphics[width=0.5\textwidth]{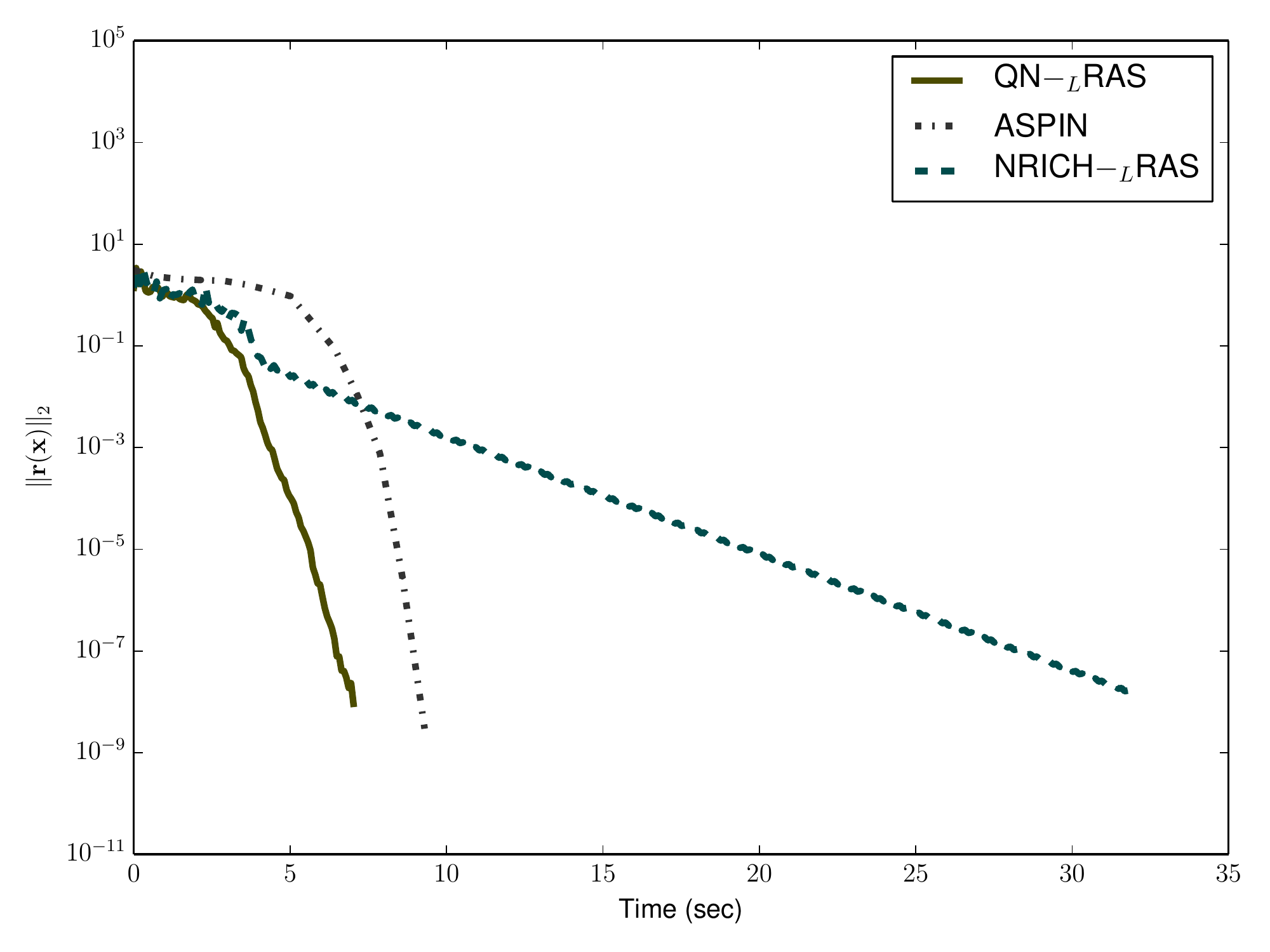}
    \caption{$\QN\lp\RAS$ , $\ASPIN$, and $\NRICH\lp\RAS$ convergence.}
    \label{fig:ex15Pre}
\end{figure}
\begin{table}[H]
  \centering
    \caption{$\QN\lp\RAS$ , $\ASPIN$, and $\NRICH\lp\RAS$ results.}
\label{tab:plaplacian:preconditioned}
\begin{tabular}{r|lllllll}
Solver & T & N. It & L. It & Func & Jac & PC & NPC \\
\hline
$\QN\lp\RAS$    & 7.02  & 92  & 0   & 410  & 185 & 185 & 185 \\
$\ASPIN$        & 9.30  & 13  & 332 & 307  & 179 & 837 & 19 \\
$\NRICH\lp\RAS$ & 32.10 & 308 & 0   & 1889 & 924 & 924 & 924 \\
\end{tabular}
\end{table}


The results for left-preconditioned methods are shown in \tabref{tab:plaplacian:preconditioned} and \figref{fig:ex15Pre}.
\ASPIN\ is competitive, taking 13 iterations and having performance characteristics similar to those of the multiplicative
\composition solver listed in \tabref{tab:plaplacian:composed}. The subdomain solvers for \ASPIN\ are set to
converge to a relative tolerance of $10^{-3}$ or 20 inner iterations.  Underresolving the local problems provides an
inadequate search direction and causes \ASPIN\ to stagnate.  The linear GMRES iteration also is converged to $10^{-3}$.

The most impressive improvements can be achieved by using a \RAS\ as a nonlinear
preconditioner.  Simple \NRICH\ acceleration does not provide much benefit compared with raw \RAS\ with respect to outer
iterations, and the preconditioner applications in the line search cause significant overhead.  However,
\QN\ using inexact \RAS\ as the left-preconditioned residual proves to be the most efficient solver for this problem,
taking 185 Newton iterations per subdomain.  Both \NRICH\ and \QN\ stagnate if the original residual is used in the line
search instead of the preconditioned one. Subdomain \QN\ methods \cite{MartinezSorSecant} have been proposed before
but built by using block approximate Jacobian inverses instead of working on the preconditioned system like \ASPIN.  As
implemented here, $\QN\lp\RAS$ and \ASPIN\ both construct left-preconditioned residuals and approximate preconditioned
Jacobian constructions; both end up being very efficient.
\section{Conclusion}
\label{sec:conclusion}
The combination of solvers using nonlinear composition, when applied carefully, may greatly improve the convergence
properties of nonlinear solvers.  Hierarchical and multilevel inner solvers allow for high arithmetic intensity and low
communication algorithms, making nonlinear composition a good option for extreme-scale nonlinear solvers.

Our experimentation, in this document and elsewhere, has shown that what works best when using nonlinear composition
varies from problem to problem.  Nonlinear composition introduces a slew of additional solver
parameters at multiple levels of the hierarchy that may be tuned for optimal performance and robustness.  A particular
problem may be amenable to a simple solver, such as \NCG, or to a combination of multilevel solvers, such as
$\FAS*(\NK\pc\MG)$. The implementation in PETSc \cite{petsc-user-ref} allows for considerable flexibility in user choice of solver compositions.

Admittedly, we have been fairly conservative in the scope of our solver combinations.  Our almost-complete restriction
to combinations of a globalized method and a decomposition method is somewhat artificial in the nonlinear case.
Without this restriction, however, the combinatorial explosion of potential methods would quickly make the scope of this
paper untenable.  Users may experiment, for their particular problem, with combinations of some number of iterations of
arbitrary combinations of nonlinear solvers, with nesting much deeper than we explore here.

The use of nonlinear preconditioning allows for solvers that are more robust to difficult nonlinear problems.  While
effective linear preconditioners applied to the Jacobian inversion problem may speed the convergence of the inner
solves, lack of the convergence of the outer Newton's method may doom the solve.  With nonlinear preconditioning, the
inner and outer treatment of the nonlinear problem allows for very rapid solution.  Nonlinear preconditioning and
\composition solvers allow for both efficiency and robustness gains, and the widespread adoption of these
techniques would reap major benefits for computational science.

\section*{Acknowledgments}
This material was based upon worked supported by the U.S. Department of Energy, Office of Science, Advanced Scientific Computing Research, under Contract DE-AC02-06CH11357. 
 We thank Jed Brown for many meaningful discussions, suggestions, and sample code.

\setcounter{tocdepth}{1}
\bibliographystyle{siam}
\bibliography{petsc.bib,petscapp.bib,composable.bib}
\newpage
{\bf Government License.} The submitted manuscript has been created by UChicago Argonne, LLC,
Operator of Argonne National Laboratory (``Argonne").
Argonne, a U.S. Department of Energy Office of Science laboratory, is
operated under Contract No. DE-AC02-06CH11357. The U.S. Government
retains for itself, and others acting on its behalf, a paid-up
nonexclusive, irrevocable worldwide license in said article to reproduce,
prepare derivative works, distribute copies to the public, and perform
publicly and display publicly, by or on behalf of the Government.

\end{document}